\documentclass[12pt]{amsart}

\usepackage{epsfig}
\usepackage{subfig}
\usepackage{amscd}
\usepackage[mathscr]{eucal}
\usepackage{amssymb}
\usepackage{amsxtra}
\usepackage{amsmath}
\usepackage[all]{xy}
\usepackage{tikz}
\usepackage{pgf}
\usetikzlibrary{calc,backgrounds,arrows,matrix,shapes}
\usepackage{verbatim}
\usepackage[bookmarks]{hyperref}

\makeatletter

\pgfdeclareshape{genuspic}{
\anchor{center}{\pgfpointorigin}
\backgroundpath{
\pgfpathmoveto{\pgfqpoint{-1cm}{0cm}}
\pgfpathcurveto %
{\pgfpoint{-0.5cm}{-.5cm}}
{\pgfpoint{0.5cm}{-.5cm}}
{\pgfpoint{1cm}{0cm}}

\pgfpathmoveto{\pgfqpoint{-0.75cm}{-0.15cm}}
\pgfpathcurveto %
{\pgfpoint{-0.25cm}{.25cm}}
{\pgfpoint{.25cm}{.25cm}}
{\pgfpoint{0.75cm}{-0.15cm}}
}
}

\pgfdeclareshape{strokegenuspic}{
\anchor{center}{\pgfpointorigin}
\backgroundpath{
\pgfpathmoveto{\pgfqpoint{-1cm}{0cm}}
\pgfpathcurveto %
{\pgfpoint{-0.5cm}{-.5cm}}
{\pgfpoint{0.5cm}{-.5cm}}
{\pgfpoint{1cm}{0cm}}

\pgfpathmoveto{\pgfqpoint{-0.75cm}{-0.15cm}}
\pgfpathcurveto %
{\pgfpoint{-0.25cm}{.25cm}}
{\pgfpoint{.25cm}{.25cm}}
{\pgfpoint{0.75cm}{-0.15cm}}
\pgfusepath{stroke}
}
}

\pgfdeclareshape{hackgenuspic}{
\anchor{center}{\pgfpointorigin}
\backgroundpath{
\pgfpathmoveto{\pgfqpoint{-0.78cm}{-.17cm}}
\pgfpathcurveto %
{\pgfpoint{-0.35cm}{-.44cm}}
{\pgfpoint{0.35cm}{-.44cm}}
{\pgfpoint{.78cm}{-0.17cm}}
\pgfpathmoveto{\pgfqpoint{-0.78cm}{-0.17cm}}
\pgfpathcurveto %
{\pgfpoint{-0.25cm}{.25cm}}
{\pgfpoint{.25cm}{.25cm}}
{\pgfpoint{0.78cm}{-0.17cm}}
\pgfsetfillcolor{white}
\pgfusepath{fill}

\pgfpathmoveto{\pgfqpoint{-1cm}{0cm}}
\pgfpathcurveto %
{\pgfpoint{-0.5cm}{-.5cm}}
{\pgfpoint{0.5cm}{-.5cm}}
{\pgfpoint{1cm}{0cm}}

\pgfpathmoveto{\pgfqpoint{-0.75cm}{-0.15cm}}
\pgfpathcurveto %
{\pgfpoint{-0.25cm}{.25cm}}
{\pgfpoint{.25cm}{.25cm}}
{\pgfpoint{0.75cm}{-0.15cm}}
\pgfusepath{stroke}
}
}

\pgfdeclareshape{fillhackgenuspic}{
\anchor{center}{\pgfpointorigin}
\backgroundpath{
\pgfpathmoveto{\pgfqpoint{-0.78cm}{-.17cm}}
\pgfpathcurveto %
{\pgfpoint{-0.35cm}{-.44cm}}
{\pgfpoint{0.35cm}{-.44cm}}
{\pgfpoint{.78cm}{-0.17cm}}
\pgfpathmoveto{\pgfqpoint{-0.78cm}{-0.17cm}}
\pgfpathcurveto %
{\pgfpoint{-0.25cm}{.25cm}}
{\pgfpoint{.25cm}{.25cm}}
{\pgfpoint{0.78cm}{-0.17cm}}
\pgfusepath{fill}

\pgfpathmoveto{\pgfqpoint{-1cm}{0cm}}
\pgfpathcurveto %
{\pgfpoint{-0.5cm}{-.5cm}}
{\pgfpoint{0.5cm}{-.5cm}}
{\pgfpoint{1cm}{0cm}}

\pgfpathmoveto{\pgfqpoint{-0.75cm}{-0.15cm}}
\pgfpathcurveto %
{\pgfpoint{-0.25cm}{.25cm}}
{\pgfpoint{.25cm}{.25cm}}
{\pgfpoint{0.75cm}{-0.15cm}}
\pgfusepath{stroke}
}
}

\makeatother

\theoremstyle{plain}

\theoremstyle{definition}

\theoremstyle{remark}

\begin{document}

\title[LG-systems and $\mathcal{D}^b(X)$ of toric Fano manifolds with small $\rho(X)$]{On Landau-Ginzburg systems and $\mathcal{D}^b(X)$ of various toric Fano manifolds with small picard group}

\address{School of Mathematical Sciences, Tel Aviv University, Ramat Aviv, Tel Aviv 6997801, Israel}
\email{yochayjerby@post.tau.ac.il}

\date{\today}

\author{Yochay Jerby}


%
%
\begin{abstract} For a toric Fano manifold $X$ denote by $Crit(X) \subset (\mathbb{C}^{\ast})^n$ the solution scheme of the Landau-Ginzburg system of equations of $X$. Examples of toric Fano manifolds with $rk(Pic(X)) \leq 3$ which admit full strongly exceptional collections of line bundles were recently found by various authors. For these examples we construct a map $E : Crit(X) \rightarrow Pic(X)$ whose image $\mathcal{E}=\left \{ E(z) \vert z \in Crit(X) \right \}$ is a full strongly exceptional collection satisfying the M-aligned property. That is, under this map, the groups $Hom(E(z),E(w))$ for $z,w \in Crit(X)$ are naturally related to the structure of the monodromy group acting on $Crit(X)$.

\end{abstract}

\maketitle

%
%

\section{Introduction and Summary of Main Results}
\label{s:intro}

\hspace{-0.6cm} Let $X$ be a smooth algebraic manifold and let $\mathcal{D}^b(X)$ be the bounded derived category of coherent sheaves on $X$, see \cite{GM,T}. A fundamental question in the study of $\mathcal{D}^b(X)$ is the question of existence of full strongly exceptional collections $\mathcal{E} = \left \{ E_1,...,E_N \right \} \subset \mathcal{D}^b(X)$. Such collections satisfy the property that the adjoint functors $$\begin{array}{ccc} R Hom_X(T, -) : \mathcal{D}^b(X) \rightarrow \mathcal{D}^b(A_\mathcal{E}) & ; & - \otimes^L_{A_{\mathcal{E}}} T : \mathcal{D}^b(A_{\mathcal{E}}) \rightarrow \mathcal{D}^b(X) \end{array} $$ are equivalences of categories where $T:=\bigoplus_{i=1}^NE_i$ and $A_{\mathcal{E}}=End(T)$ is the corresponding endomorphism ring. The first example of such a collection is $$\mathcal{E} = \left \{ \mathcal{O}, \mathcal{O}(1),...,\mathcal{O}(s) \right \} \subset Pic(\mathbb{P}^s)$$ found by Beilinson in \cite{B}. When $X$ is a toric manifold one further asks the more refined question of whether $\mathcal{D}^b(X)$ admits an exceptional collection whose elements are line bundles $\mathcal{E} \subset Pic(X)$, rather than general elements of $\mathcal{D}^b(X)$?

\hspace{-0.6cm} Let $X$ be a $s$-dimensional toric Fano manifold given by a Fano polytope $\Delta$ and and let $\Delta^{\circ}$ be the polar polytope of $\Delta$. Let $f_X = \sum_{n \in \Delta^{\circ} \cap \mathbb{Z}^n} z^n \in \mathbb{C}[z_1^{\pm},...,z_s^{\pm}]$ be the Landau-Ginzburg potential associated to $X$, see \cite{Ba,FOOO,OT}. Recall that the Landau-Ginzburg system of equations is given by $$ z_i \frac{\partial}{\partial z_i } f_X(z_1,...,z_s)=0 \hspace{0.5cm} \textrm{ for } \hspace{0.25cm} i=1,...,s $$ and denote by $Crit(X) \subset (\mathbb{C}^{\ast})^s$ the corresponding solution scheme.

\hspace{-0.6cm} In \cite{J} we defined an exceptional map to be a map of the form 
$E : Crit(X) \rightarrow Pic(X)$ such that $\mathcal{E}_E(X) := E(Crit(X)) \subset Pic(X)$ is a full strongly exceptional collection of line bundles and introduced natural exceptional maps for the five Del-Pezzo surfaces and various three dimensional toric Fano manifolds. Our main observation in \cite{J} was that, in the considered cases, the defined exceptional maps exhibit a non-trivial property to which we refer as the M-aligned property. In short, M-aligned means that the algebraic structure of the spaces $Hom(E(z),E(w))$ for $z,w \in Crit(X)$ could be further realted to the structure of the geometric monodromy group acting on $Crit(X)$.

\hspace{-0.6cm} In general, not much is known on the question of which toric Fano manifolds admit full strongly exceptional collections of line bundles in $Pic(X)$. And, with full generality, the question seems out of reach with current techniques. However, full strongly exceptional collections for specific classes of toric Fano manifolds with small $\rho(X) = rank(Pic(X))$ were recently found by various authors, see \cite{CMR3,CMR4,DLM}. Our aim in this work is to introduce exceptional maps $E : Crit(X) \rightarrow Pic(X) $ for these classes and study the $M$-aligned property for these maps.

\hspace{-0.6cm} Defining a map of the form $E: Crit(X) \rightarrow Pic(X)$ requires the association of integral invariants to the elements of the solution scheme $Crit(X)$. The Landau-Ginzburg potential $f_X$ is an element of the space $$ L(\Delta^{\circ}) : = \left \{ \sum_{n \in \Delta^{\circ} \cap \mathbb{Z}^n} u_n z^n \vert u_n \in \mathbb{C}^{\ast} \right \} \subset \mathbb{C}[z_1^{\pm},...,z_n^{\pm}]$$ Similarly, for $f_u(z) = \sum_n u_n z^n \in L(\Delta^{\circ})$ define the corresponding Landau-Ginzburg system of equations $$ z_i \frac{\partial}{\partial z_i } f_u(z_1,...,z_n) =0 \hspace{0.5cm} \textrm{ for } \hspace{0.25cm} i=1,...,n $$ and denote by $Crit(X; f_u ) \subset (\mathbb{C}^{\ast})^n$ the solution scheme. Our first main observation, is that the arguments of the solution scheme $Crit(X; f_u)$ converge to roots of unity as $Log \vert u \vert \rightarrow \pm \infty$. We view these collections of roots of unity as giving an "asymptotic" generalization of the roots of unity arising directly in the case of projective space.

\hspace{-0.6cm} Denote by $ \Delta(k)$ the set of $k$-dimensional facets of $\Delta$. Note that $ \Delta( n-1) \simeq \Delta^{\circ}(0)$ and denote by $F_n \in \Delta(n-1)$ the facet corresponding to $n \in \Delta^{\circ} (0)$. Denote by $V_X(F_n)$ the $T$-orbit closure corresponding to the facet $F_n$. In particular, for a solution $z \in Crit(X ; f_u) $ define the $T$-invariant $\mathbb{R}$-divisor $$ D(z) := \sum_{n \in \Delta^{\circ} \cap \mathbb{Z}^n} Arg(z^n) \cdot V_X(F_n) \in Div_T(X) \otimes \mathbb{R} $$ In the cases considered, $Pic(X)$ admits a natural set of generators. In particular, denote by $\left [D(z) \right ] \in Pic(X)$ the integral part of $D(z)$ via the map $Div_T(X) \otimes \mathbb{R} \rightarrow Pic(X) \otimes \mathbb{R}$. In the first part of the paper we show:

\bigskip

\hspace{-0.6cm} \bf Theorem A: \rm Let $X$ be a toric Fano manifold of one of the following classes:

\bigskip

(a) $X = \mathbb{P}(\mathcal{O}_{\mathbb{P}^s}\oplus \left (\oplus_{i=1}^r\mathcal{O}_{\mathbb{P}^s}(a_i) \right ) )$, a projective bundle with $\sum_{i=1}^r a_i \leq s$.

\bigskip

(b) $X = Bl_B (\mathbb{P}^{n-r} \times \mathbb{P}^r)$, where $B=\mathbb{P}^{n-r-1} \times \mathbb{P}^{r-1}$.

\bigskip

(c) $X=Bl_B(\mathbb{P}(\mathcal{O}_{\mathbb{P}^{n-1}} \oplus \mathcal{O}_{\mathbb{P}^{n-1}}(b)))$, where $B=\mathbb{P}^{n-2}$ and $b<n-1$.

\bigskip

\hspace{-0.6cm} Then $X$ admits an exceptional map $E : Crit(X ; f_t) \rightarrow Pic(X)$ given by $z \mapsto \left [ D(z) \right ]$ for $0<<t$ big enough.

\bigskip

\hspace{-0.6cm} Note that for (a) one has $\rho(X)=rk(Pic(X)) =2$ while for (b) and (c) one has $\rho(X) =3$. In fact, due to Kleindschmit's classification theorem any toric manifold with $\rho(X)=2$ is of class (a), see \cite{K}. The exceptional collections arising for (a) coincide with those found by Costa and Miro-Roig in \cite{CMR3}. The exceptional collections arsing for classes (b) and (c) coincide with those found by Costa and Miro-Roig in \cite{CMR4}
and by Lason, Michalek and Dey in \cite{DLM}, respectively. Both of which are special cases of a more general construction for a wider class of toric manifolds with $\rho(X) =3$ due to \cite{LM}. Let us illustrate the above with the following example:

\bigskip

\hspace{-0.6cm} \bf Example \rm (projective line): For $X= \mathbb{P}^s$ the Landau-Ginzburg potential is given by $f(z_1,...,z_s)=z_1+...+z_s +\frac{1}{z_1 \cdot ... \cdot z_s}$ and the corresponding system
of equations is $$ z_i \frac{\partial}{\partial z_i } f_X(z_1,...,z_s)=z_i - \frac{1}{z_1 \cdot ... \cdot z_s} =0 \hspace{0.5cm} \textrm{ for } \hspace{0.25cm} i=1,...,s $$ The solution scheme $Crit(\mathbb{P}^s) \subset (\mathbb{C}^{\ast})^s$ is given by
$z_k= ( e^{\frac{2 \pi ki}{s+1}},...,e^{\frac{2 \pi k i}{s+1}}) $ for $k=0,...,s$, via direct computation. We have $Div_T(X) = \bigoplus_{i=1}^r V_X(F_i) \cdot \mathbb{Z}$ where $V_X(F_i)$ is the projective hyperplane defined by $z_i=0$. The exceptional map is hence given by $$E(z_k) = \left [\frac{k}{s+1} \cdot V_X(F_1) + ... + \frac{k}{s+1} \cdot V_X(F_s) \right ] = k \cdot H$$ where $H \in Pic(X)$ is the positive generator.
\bigskip

\hspace{-0.6cm} A full strongly exceptional collection of line bundles $\mathcal{E} \subset Pic(X) $ on a toric manifold $X$ is associated with a quiver $Q_{\mathcal{E}}$ whose vertices are the elements of $\mathcal{E}$ and edges connecting two elements $E_1,E_2 \in \mathcal{E}$ correspond to the $T$-invariant divisors appearing in the natural splitting $$ Hom(E_1,E_2) \simeq \bigoplus_{D \in Div_T(X)} Hom_D(E_1,E_2)$$ The main feature of the exceptional maps, defined in \cite{J}, is the $M$-aligned property, which relates the quivers $Q_{\mathcal{E}}$ to monodromies acting on the solution scheme $Crit(X; f_t)$. The second part of the paper is devoted for the description of the $M$-aligned for the exceptional maps introduced for the classes (a)-(c) in Theorem A. For instance, for a projective bundle $X = \mathbb{P}( \mathcal{O} \oplus ( \oplus_{i=1}^r \mathcal{O}(a_i)))$ the corresponding potential function is given by $$ f_t(z,w) :=1+ \sum_{i=1}^s z_i + \sum_{i=1}^r w_i + e^{t} \cdot \frac{w_1^{a_1} \cdot ... \cdot w_r^{a_r}}{z_1 \cdot ... \cdot z_s} + \frac{1}{w_1 \cdot ... \cdot w_r}$$ for a $T$-divisor $$D:= \sum_{i=1}^s n_i \cdot V_X(e_i) + \sum_{i=s+1}^{s+r} m_i \cdot V_X(e_i)+ n_0 \cdot V_X \left ( \sum_{i=s+1}^r a_i e_i - \sum_{i=1}^s e_i \right) + m_0 \cdot V_X \left (-\sum_{i=s+1}^r e_i \right ) $$ we associate following loop of Laurent polynomials in $L(\Delta^{\circ})$: $$ \gamma_D^{\theta}(z,w) :=1+ \sum_{i=1}^s e^{2 \pi i n_i \theta} z_i + \sum_{i=1}^r e^{2 \pi i m_i \theta} w_i + e^{t+2 \pi i n_0 \theta} \cdot \frac{w_1^{a_1} \cdot ... \cdot w_r^{a_r}}{z_1 \cdot ... \cdot z_s} +e^{2 \pi i m_0 \theta} \frac{1}{w_1 \cdot ... \cdot w_r}$$ for $ \theta \in [0,1]$. Note that $\gamma^0_D = \gamma^1_D$. In particular, following the solution schemes $Crit(X ; \gamma_D^{\theta})$ along $\theta \in [0,1]$ one gets a map $M_D : Crit(X ; f_t) \rightarrow Crit(X ; f_t)$, to which we refer as the monodromy map corresponding to the divisor $D$. We show:

\bigskip

\hspace{-0.6cm} \bf Theorem B \rm ($M$-aligned property for projective bundles): Let $z_1,z_2 \in Crit(X ; f_t)$ be two elements and $D \in Div_T(X)$ be a $T$-divisor. Then $$ Hom_D(E(z_1),E(z_2)) \neq 0 \Rightarrow M_D(z_1) = z_2$$ where $E : Crit(X ; f_t) \rightarrow Pic(X)$ is the exceptional map of Theorem A.

\bigskip

\hspace{-0.6cm} We refer to this property as $M$-algined to indicate that the $Hom$-spaces between elements of the collection $\mathcal{E}$ are algined with the corresponding monodromy actions on $Crit(X ; f_t)$. In section 4 we prove Theorem B and further consider the $M$-aligned property for the classes (b) and (c) of Theorem A, as well.

\bigskip

\hspace{-0.6cm} The rest of the work is organized as follows: In section 2 we recall relevant facts on toric Fano manifolds and their derived categories of coherent sheaves and recall the definition of the exceptional collections for the classes (a)-(c). In section 3 we study the Landau-Ginzburg system for classes (a)-(c) and define the exceptional map. In section 4 we show the M-aligned property. In section 5 we discuss concluding remarks and relations to further topics of mirror symmetry.

\section{Relevant Facts on Toric Fano Manifolds}
\label{s:Rfotfm}

\hspace{-0.6cm} Let $N \simeq \mathbb{Z}^n$ be a lattice and let $M = N^{\vee}=Hom(N, \mathbb{Z})$ be the dual lattice. Denote by $N_{\mathbb{R}} = N \otimes \mathbb{R}$ and $M_{\mathbb{R}}=M \otimes \mathbb{R}$ the corresponding vector spaces. Let $ \Delta \subset M_{\mathbb{R}}$ be an integral polytope and let $$ \Delta^{\circ} = \left \{ n \mid (m,n) \geq -1 \textrm{ for every } m \in \Delta \right \} \subset N_{\mathbb{R}}$$ be the \emph{polar} polytope of $\Delta$. The polytope $\Delta \subset M_{\mathbb{R}}$ is said to be \emph{reflexive} if $0 \in \Delta$ and $\Delta^{\circ} \subset N_{\mathbb{R}}$ is integral. A reflexive polytope $\Delta$ is said to be \emph{Fano} if every facet of $\Delta^{\circ}$ is the convex hall of a basis of $M$.

\hspace{-0.6cm} To an integral polytope $\Delta \subset M_{\mathbb{R}}$ associate the space $$ L(\Delta) = \bigoplus_{m \in \Delta \cap M} \mathbb{C} m $$ of Laurent polynomials whose Newton polytope is $\Delta$. Denote by $i_{\Delta} : (\mathbb{C}^{\ast})^n \rightarrow \mathbb{P}(L(\Delta)^{\vee})$ the embedding given by $ z \mapsto [z^m \mid m \in \Delta \cap M] $. The \emph{toric variety} $X_{\Delta} \subset \mathbb{P}(L(\Delta)^{\vee})$ corresponding to the polytope $\Delta \subset M_{\mathbb{R}}$ is defined
to be the compactification of the image $i_{\Delta}((\mathbb{C}^{\ast})^n) \subset \mathbb{P}(L(\Delta)^{\vee})$. A toric variety $X_{\Delta}$ is said to be Fano if its anticanonical class $-K_X$ is Cartier and ample. In \cite{Ba2} Batyrev shows that $X_{\Delta}$ is a Fano variety if $\Delta$ is reflexive and, in this case, the embedding
$i_{\Delta}$ is the anti-canonical embedding. The Fano variety $X_{\Delta}$ is smooth if and only if $\Delta^{\circ}$ is a Fano polytope.

\hspace{-0.6cm} Denote by $\Delta(k)$ the set of $k$-dimensional faces of $\Delta$ and denote by $ V_X(F) \subset X$ the orbit closure of the orbit corresponding to the facet $F \in \Delta(k)$ in $X$, see \cite{F,O}. In particular, consider the group of toric divisors $$ Div_T(X) := \bigoplus_{F \in \Delta(n-1)} \mathbb{Z} \cdot V_X(F)$$ Assuming $X$ is smooth the group $Pic(X)$ is described in terms of the short exact sequence $$ 0 \rightarrow M \rightarrow
Div_T(X) \rightarrow Pic(X) \rightarrow 0 $$ where the map on the left hand side is given by $ m \rightarrow \sum_F \left < m, n_F \right > \cdot V_X(F) $ where $n_F \in \mathbb{N}_{\mathbb{R}}$ is the unit normal to the hyperplane spanned by the facet $F \in \Delta(n-1)$. In particular, note that $$ \rho(X)= rank \left ( Pic (X) \right ) = \vert \Delta(n-1) \vert -n $$ Moreover, when $\Delta$ is reflexive one has $\Delta^{\circ}(0)=\left \{ n_F \vert F \in \Delta(n-1) \right \} \subset N_{\mathbb{R}}$. We thus sometimes denote $V_X(n_F)$ for the $T$-invariant divisor $V_X(F)$. We denote by $Div_T^+(X) $ the semi-group of all toric divisors $ \sum_{F} m_F \cdot V_X(F)$ with $0 \leq m_F$ for any $F \in \Delta(n-1)$.

\hspace{-0.6cm} Let $X$ be a smooth projective variety and let $\mathcal{D}^b(X)$ be the derived category of bounded complexes of
coherent sheaves of $\mathcal{O}_X$-modules, see \cite{GM,T}. For a finite dimensional algebra $A$ denote by $\mathcal{D}^b(A)$ the derived category of bounded complexes of finite dimensional right modules over $A$. Given an object $T \in \mathcal{D}^b(X)$ denote by $A_T=Hom(T,T)$ the corresponding endomorphism algebra.

\bigskip

\hspace{-0.6cm} \bf Definition 2.1: \rm An object $T \in \mathcal{D}^b(X)$ is called a \emph{tilting object} if the corresponding adjoint functors $$\begin{array}{ccc} R Hom_X(T, -) : \mathcal{D}^b(X) \rightarrow \mathcal{D}^b(A_T) & ; & - \otimes^L_{A_T} T : \mathcal{D}^b(A_T) \rightarrow \mathcal{D}^b(X) \end{array} $$ are equivalences of categories. A locally free tilting object is called a tilting bundle.

\bigskip

\hspace{-0.6cm} An object $ E \in \mathcal{D}^b(X)$ is said to be \emph{exceptional} if $Hom(E,E)=\mathbb{C}$ and $Ext^i(E,E)=0$ for $0<i$. We have:

\bigskip

\hspace{-0.6cm} \bf Definition 2.2: \rm An ordered collection $ \mathcal{E} = \left \{ E_1,...,E_N \right \} \subset \mathcal{D}^b(X)$ is said to be an \emph{exceptional collection} if each $E_j$ is exceptional and $Ext^i(E_k,E_j) =0$ for $j<k \textrm{ and } 0 \leq i $. An exceptional collection is said to be \emph{strongly exceptional} if also $Ext^i(E_j,E_k)=0$ for $j \leq k$ and $0<i$. A strongly exceptional collection is called \emph{full} if its elements generate $\mathcal{D}^b(X)$ as a triangulated category.

\bigskip

\hspace{-0.6cm} The importance of full strongly exceptional collections in tilting theory is due to the following properties, see \cite{Bo,K}:

\bigskip

- If $ \mathcal{E}$ is a full strongly exceptional collection then $T = \bigoplus_{i=1}^N E_i$ is a tilting object.

\bigskip

- If $T = \bigoplus_{i=1}^N E_i$ is a tilting object and $ \mathcal{E} \subset Pic(X) $ then $\mathcal{E}$ can be ordered as a full

\hspace{0.2cm} strongly exceptional collection of line bundles.

\bigskip

\hspace{-0.6cm} Let us note the following examples of toric Fano manifolds, which would be considered in the continuation:

\bigskip

\hspace{-0.6cm} \bf Example 2.3 \rm (Projective bundles): By a result of Kleinschmidt's \cite{Kl} the class of toric manifolds with $\rho(X) = rk(Pic(X))=2$ consists of the projective bundles $$ X_a=\mathbb{P} \left (\mathcal{O}_{\mathbb{P}^s} \oplus \bigoplus_{i=1}^r \mathcal{O}_{\mathbb{P}^s}(a_i) \right ) \hspace{0.5cm} \textrm{ with } \hspace{0.25cm} 0 \leq a_1 \leq... \leq a_r$$ see also \cite{CoLS}. Set $a_0 =0$. Consider the lattice $N = \mathbb{Z}^{s+r}$ and let $v_1,...,v_s$ be the standard basis elements of $\mathbb{Z}^s$ and $e_1,...,e_r$ be the standard basis elements of $\mathbb{Z}^r$. Set $v_0 =- \sum_{i=1}^s u_i+ \sum_{i=1}^s a_i e_i $ and $ e_0 =- \sum_{i=1}^r e_i$. Let $\Delta_a^{\circ} \subset N_{\mathbb{R}}$ be the polytope whose vertex set is given by $$\Delta_a^{\circ}(0) = \left \{ v_0,...,v_s,e_0,...,e_r \right \}$$ It is straightforward to verify that $\Delta^{\circ}_a$ is a Fano polytope if and only if $\sum_{i=1}^r a_i \leq s$. In particular, in this case $X_a \simeq X_{\Delta_a}$, see \cite{CoLS}. One has $$ Pic(X_a) = \xi \cdot \mathbb{Z} \oplus \pi^{\ast}H \mathbb{Z}$$ where $\xi$ is the class of the tautological bundle and $\pi^{\ast}H $ is the pullback of the generator $H$ of $Pic(\mathbb{P}^s) \simeq H \cdot \mathbb{Z}$ under the projection $\pi : X_a \rightarrow \mathbb{P}^s$. Note that the following holds $$ \begin{array}{ccc} [V_X(v_i)] = \pi^{\ast}H & ; & [V_X(e_j)]=\xi - a_i \cdot \pi^{\ast} H \end{array}$$ for $0 \leq i \leq s $ and $ 0 \leq j \leq r$. In \cite{CMR3} Costa and Mir$\acute{\textrm{o}}$-Roig show that the collection of line bundles $\mathcal{E} = \left \{E_{kl} \right \}_{k=0,l=0}^{s,r} \subset Pic(X)$ where $$ E_{kl} := k \cdot \pi^{\ast} H + l \cdot \xi \hspace{0.5cm} \textrm{ for } \hspace{0.25cm} 0 \leq k \leq s \hspace{0.1cm} , \hspace{0.1cm} 0 \leq l \leq r $$ admits the structure of a full strongly exceptional collection.

\bigskip

\hspace{-0.6cm} \bf Example 2.4: \rm For $r \leq n $ let $ X_{n,r}=Bl_B(\mathbb{P}^{n-r} \times \mathbb{P}^r)$ be the blow up of $\mathbb{P}^{n-r} \times \mathbb{P}^r$ along the multi-linear codimension two subspace $B = \mathbb{P}^{n-r-1} \times \mathbb{P}^{r-1}$. Let $e_1,...,e_n$ be the standard basis of $M_{\mathbb{R}}$. The vertices of the polar polytope are given by $ \Delta^{\circ}(0) = \left \{ e_1,...,e_n, v_1 , v_2, v_3 \right \}$ where $$ \begin{array}{ccccc} v_1 = -\sum_{i=1}^{n-r} e_i & ; & v_2 = -\sum_{i=n-r+1}^n e_i & ; & v_3 = -\sum_{i=1}^n e_i \end{array} $$ In particular, $$ Pic(X_{n,r}) = U \cdot \mathbb{Z} \oplus V \cdot \mathbb{Z} \oplus E \cdot \mathbb{Z}$$ with $$\begin{array}{ccccc} U:=[V_X(e_i)]=[ V_X(v_1) + V_X(v_3) ] & ; & V:=[V_X(e_j)]=[V_X(v_2) +V_X(v_3)] & ; & E: = [V_X(v_3)] \end{array} $$ where $1 \leq i \leq n-r$ and $n-r+1 \leq j \leq n$. Set $$ E_{kl}:= kU+lV \hspace{0.5cm} \textrm{ for } \hspace{0.25cm} 0 \leq k \leq n-r \textrm{ , } 0 \leq l \leq r$$ and $$ F_{m,n} :=mU +n V - E \hspace{0.5cm} \textrm{ for } \hspace{0.25cm} 1 \leq m \leq n-r \textrm{ , } 1 \leq n \leq r$$ In \cite{CMR4} Costa and Miro-Roig show that the following collection $$ \mathcal{E}_{n,r} = \left \{E_{kl} \right \}_{k,l=0}^{n-r,r} \cup \left \{ F_{m,n} \right \}_{m,n=1}^{n-r,r} \subset Pic(X_{n,r})$$ admits the structure of a full strongly exceptional collection.

\bigskip

\hspace{-0.6cm} \bf Example 2.5: \rm For $b < n-1$ let $X_{n,b} = Bl_B ( \mathbb{P}( \mathcal{O}_{\mathbb{P}^{n-1}} \oplus \mathcal{O}_{\mathbb{P}^{n-1}}(b)))$ the blow-up of the projective bundle $ \mathbb{P}( \mathcal{O}_{\mathbb{P}^{n-1}} \oplus \mathcal{O}_{\mathbb{P}^{n-1}}(b))$ along $B \simeq \mathbb{P}^{n-2}$. The vertices of the polar polytope are given by $ \Delta^{\circ}(0) = \left \{ e_1,...,e_n, u_1 , u_2, u_3 \right \}$ where $$ \begin{array}{ccccc} u_1 = -e_n & ; & u_2 = - \sum_{i=1}^{n-1} e_i- be_n & ; & u_3 = -\sum_{i=1}^{n-1} e_i - (b+1)e_n \end{array} $$ We have $$ Pic(X) \simeq V \cdot \mathbb{Z} \oplus Y \cdot \mathbb{Z} \oplus T \cdot \mathbb{Z}$$ where $V=[V_X(e_i)]$ for $i=1,...,n-1$, $Y=[V_X(u_3)]$ and $T=[V_X(u_1)]$. For $0 \leq k \leq n-1$, $ 0 \leq l \leq 1 $ and $1 \leq m \leq n-1$ set $$ \begin{array}{ccc} E_{kl} :=k V + l(Y+T+bV) & ; & F_m := (m+b) V + T \end{array} $$ Michalek, Lason and Dey show in \cite{DLM} that $$ \mathcal{E}_{n,b} := \left \{E_{kl} \right \}_{k=0,l=0}^{n-1,1} \cup \left \{F_m \right \}_{m=1}^{n-1} \subset Pic(X) $$ admits the structure of a full strongly exceptional collection.

\bigskip

\hspace{-0.6cm} Let us conclude this section by recalling a few facts about the Frobenius toric map and Frobenius splitting:

\bigskip

\hspace{-0.6cm} \bf Remark 2.6 \rm (Exceptional collections and Frobenius splitting): Let $ l \in \mathbb{N}$ be an integer and let $F_l : (\mathbb{C}^{\ast})^s \rightarrow (\mathbb{C}^{\ast})^s$ be the $l$-th power map given by $(x_1,...,x_s) \mapsto (x_1^l,...,x_s^l)$. The extension of this map to the whole toric manifold $X$, which we also denote by $F_l : X \rightarrow X $, is called the $l$-th toric Frobenius mapping. In \cite{T} Thomoson described the following formula for the Frobenius splitting of the trivial bundle $$ (F_l)_{\ast} \mathcal{O} = \bigoplus_{D \in Pic(X)} \mathcal{O}(D)^{m(D)} $$ where $m(D)$ is the number of points in the cube $Div_T(X)/l \cdot Div_T(X)$ representing the class $-lD \in Pic(X)$, see also \cite{A}. Let $\widetilde{F}_l : Div_T(X)/ l \cdot Div_T(X) \rightarrow Pic(X)$ be the function given by $ D \mapsto [D/l]$. In particular, in the limit as $l \rightarrow \infty$ we get a map $\widetilde{F} : \mathbb{T}^{r} \rightarrow Pic(X)$, where $r = rk(Div_T(X))$, see \cite{B}. Let $\mathcal{B}_X \subset Pic(X)$ be the image of the map $\widetilde{F}$. The exceptional collection of Example 2.3-4 satisfy $\mathcal{E}_X = \mathcal{B}_X$ while the exceptional collection of Example 4.5 satisfies $\mathcal{E}_X \subset \mathcal{B}_X$.

\hspace{-0.6cm} Let $\mathcal{B}_X \subset Pic(X)$ be the image set of the map $\widetilde{F}$. In general, $\chi(X) \leq \vert \mathcal{B}_X \vert$. Hence, the number of elements in $\mathcal{B}_X$ exceeds the expected number of elements of an exceptional collection. On the one hand, many results from recent years seem to indicate an ilusive relationship between the set $\mathcal{B}_X$ and full strongly exceptional collections of line bundles, see \cite{Bo3,CMR2,CMR4,LM}. On the other hand, Michalek, Lason and Dey found examples of toric Fano manifolds for which the set $\mathcal{B}_X$ admits no subset which is a full strongly exceptional collection, see \cite{LM}. An explanation for the distinction between the two cases is yet unclear.

\section{Variations of the LG-system and exceptional maps}
\label{s:LGM}

\hspace{-0.6cm} Let $X$ be a $n$-dimensional toric Fano manifold given by a Fano polytope $\Delta \subset M_{\mathbb{R}}$ and let $\Delta^{\circ} \subset N_{\mathbb{R}}$ be the corresponding polar polytope. Set $$L(\Delta^{\circ}):= \left \{ \sum_{n \in \Delta^{\circ} \cap \mathbb{Z}^n} u_n z^n \vert u_n \in \mathbb{C}^{\ast} \right \} \subset \mathbb{C}[z_1^{\pm},...,z_n^{\pm}]$$ We refer to $$ z_i \frac{\partial}{\partial z_i } f_u(z_1,...,z_n)=0 \hspace{0.5cm} \textrm{ for } \hspace{0.25cm} i=1,...,n $$ as the LG-system of equations associated to an element $ f_u(z) = \sum_{n \in \Delta^{\circ} \cap \mathbb{Z}^n} u_n z^n $ and denote by $Crit(X; f_u) \subset (\mathbb{C}^{\ast})^n$ the corresponding solution scheme. We refer to the element $f_X(z) = \sum_{n \in \Delta^{\circ} \cap \mathbb{Z}^n} z^n$ as the LG-potential of $X$.

\bigskip

\hspace{-0.6cm} \bf Definition 3.1: \rm A map $E : Crit(X; f_u) \rightarrow Pic(X)$ is an exceptional map if the image $E(Crit(X; f_u)) \subset Pic(X)$ is a full strongly exceptional collection.

\bigskip

\hspace{-0.6cm} Our main observation in this section is that roots of unity arise when considering the asymptotics of $Crit(X ; f_u)$ as $Log \vert u \vert \rightarrow \pm \infty$. We view these roots of unity as a generalization of the roots of unity arising in $Crit(X)$ in the case of projective space $X= \mathbb{P}^n$.

\hspace{-0.6cm} Let $Arg: (\mathbb{C}^{\ast})^n \rightarrow \mathbb{T}^n$ be the map given by $ (r_1 e^{2 \pi i \theta_1},..., r_n e^{2 \pi i \theta_n}) \mapsto (\theta_1 ,...,\theta_n )$. For any $z \in Crit(X ; f_u) $ define the $\mathbb{R}$-divisor $$ D_u(z) := \sum_{n \in \Delta^{\circ}(0)} Arg(z^n) \cdot V_X(F_n) \in Div_T(X) \otimes \mathbb{R}$$ Assume $L_1,...,L_{\rho} \in Pic(X)$ is a given basis. Let $D = \sum_n a_n \cdot V_X(F_n) \in Div_T(X) \otimes \mathbb{R}$ be an $\mathbb{R}$-divisor and denote by $ [D] := \sum_{i=1}^{\rho} b_i \cdot L_i$ the corresponding element of $Pic(X) \otimes \mathbb{R}$. We denote by $$ [D]_{\mathbb{Z}} := \sum_{i=1}^{\rho} [b_i]_{\mathbb{Z}} \cdot L_i \in Pic(X)$$ the line bundle obtained by replacing the coefficients of $[D]$ with their integer part. Define the map $E_u : Crit(X ; f_u) \rightarrow Pic(X)$ given by $z \mapsto [D_u(z)]_{\mathbb{Z}} \in Pic(X)$. In the continuation of this section we show how maps of the form $E_u$ give rise to exceptional maps for the toric Fano manifolds of examples 2.3-2.5 ((a)-(c) of the introduction).
\bigskip

\hspace{-0.6cm} \bf Remark 3.1 \rm (Geometric viewpoint): Let $\left \{f_u \right \}_{u \in \mathbb{C}^{\ast}} \subset L(\Delta^{\circ})$ be a $1$-parametric family of Laurent polynomials. Consider the Riemann surface $$ C: = \left \{ \left ( z ,u \right ) \vert z \in Crit(X; f_u) \right \} \subset (\mathbb{C}^{\ast})^n \times \mathbb{C}^{\ast}$$ Denote by $\pi : C \rightarrow \mathbb{C}^{\ast}$ the projection on the second factor which expresses $C$ as an algebraic fibration over $\mathbb{C}^{\ast}$ of rank $N=\chi(X)$. Denote by $C(u) = \pi^{-1}(u)$ for $u \in \mathbb{C}^{\ast}$. A graphic illustration of $C$ together with the curves $C(t)$ for $0 \leq t$ for the Hirzebruch surface $X= \mathbb{P}(\mathcal{O}_{\mathbb{P}^1} \oplus \mathcal{O}_{\mathbb{P}^1}(1))$ is as follows:
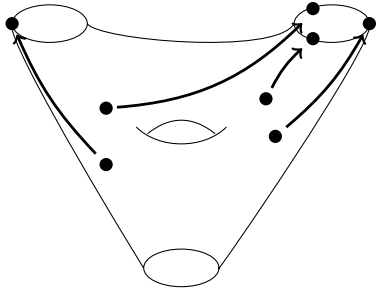
\begin{center}
$$ \begin{tikzpicture}[scale=.25]

\node [genuspic, draw, scale=0.6, fill=white] at (0,-2.5) {};

\draw (-7,3) ellipse (2 cm and 1 cm)
(-5,3) .. controls +(-60:1) and +(-120:2) .. (6,3)

(8,3) ellipse (2cm and 1cm)
(10,3) .. controls +(-60:1) and +(-120:1) .. (2,-10)
(-0,-10) ellipse (2cm and 1cm)
(-2,-10).. controls +(-60:1) and +(-120:1) .. (-9,3);

\node (a) at (5,-3) {};
\node (b) at (-4,-4.5){};
\node (b2) at (7,3.6) {};
\node (c) at (-4,-1.5){} ;
\node (d) at (4.5,-1){} ;
\fill[black] (a) circle (10pt);
\fill[black] (b) circle (10pt);
\fill[black] (c) circle (10pt);
\fill[black] (d) circle (10pt);

\node (a1) at (10,3) {};
\node (b1) at (7,3.8){};
\node (c1) at (7,2.2){} ;
\node (d1) at (-9,3){} ;
\fill[black] (a1) circle (10pt);
\fill[black] (b1) circle (10pt);
\fill[black] (c1) circle (10pt);
\fill[black] (d1) circle (10pt);

\path[every node/.style={font=\sffamily\small}]
(a) edge [->,bend right=10,line width=1.05pt] (a1)

(b) edge [->,bend left=10,line width=1.05pt](d1)

(c) edge [->,bend right=20,line width=1.05pt](b2)
(d) edge [->,bend left=10,line width=1.05pt](c1);
\end{tikzpicture}
$$
\end{center}

\hspace{-0.6cm} An amusing analogy can be drawn between the resulting dynamics and the cue game of "pool". Indeed, consider the Riemann surface $C$ as a "pool table", the cusps of the surface as the "pockets", and the set $C(0) \simeq Crit(X)$ as an initial set of "balls". In this analogy the dynamics of $C(t)$ describes the path in which the balls approach the various "pockets" of the table as $ t \rightarrow \pm \infty $.
\bigskip

\subsection{Exceptional map for projective bundles} For the projective bundle $X_a$ the Landau-Ginzburg potential is given by $$ f(z,w)=1+\sum_{i=1}^s z_i+\sum_{i=1}^r w_i + \frac{w_1^{a_1} \cdot ... \cdot w_r^{a_r}}{z_1 \cdot ... \cdot z_s } + \frac{1}{w_1 \cdot ... \cdot w_r } \in L(\Delta_a^{\circ})$$ We consider the $1$-parametric family of Laurent polynomials $$ f_u(z,w):=1+\sum_{i=1}^s z_i+\sum_{i=1}^{r} w_i + e^{u} \cdot \frac{w_1^{a_1} \cdot ... \cdot w_r^{a_r}}{z_1 \cdot ... \cdot z_s } + \frac{1}{w_1 \cdot ... \cdot w_r } \in L(\Delta_a^{\circ})$$ for $u \in \mathbb{C}$.

\hspace{-0.6cm} Let $A(V):=Arg(V) \subset \mathbb{T}^n$ be the image of the algebraic subvariety $V \subset (\mathbb{C}^{\ast})^n $ under the argument map. Such sets are known as \emph{co-amoebas}, see \cite{PT}. For $1 \leq i \leq s$ and $1 \leq j \leq r$ consider the following sub-varieties of $(\mathbb{C}^{\ast})^{s+r}$: $$ \begin{array}{ccc} V_i^u = \left \{ z_i - e^u \frac{ \prod_{i=1}^r w_i^{a_i} }{\prod_{i=1}^s z_i }= 0 \right \} & ; & W_j^u = \left \{ w_i +a_i e^u \frac{ \prod_{i=1}^r w_i^{a_i} }{\prod_{i=1}^s z_i}-\frac{1}{ \prod_{i=1}^r w_i} = 0 \right \} \end{array} $$ For the co-amoeba one has $$ A(Crit(X ; f_u) ) \subset \left (\bigcap_{i=1}^s A(V_i^u) \right ) \cap \left ( \bigcap_{i=1}^r A(W_i^u) \right ) \subset \mathbb{T}^{s+r} $$ Let $( \theta_1,...,\theta_s, \delta_1,...,\delta_r)$ be coordinates on $\mathbb{T}^{s+r}$. We have, via straight-forward computation:

\bigskip

\hspace{-0.6cm} \bf Lemma 3.1.1: \rm For $1 \leq i \leq s$ and $1 \leq j \leq r$:

\bigskip

(1) $lim_{t \rightarrow -\infty} A(V_i^t)= \left \{ \theta_i + \sum_{i=1}^s \theta_i - \sum_{j=1}^r a_j \delta_j =0 \right \} \subset \mathbb{T}^{s+r}$

\bigskip

(2) $ lim_{t \rightarrow -\infty} A(W_j^t)= \left \{ \delta_j + \sum_{j=1}^r \delta_j=0 \right \} \subset \mathbb{T}^{s+r}$

\bigskip

\hspace{-0.6cm} Let $ \Theta : (\mathbb{C}^{\ast})^{s+r} \rightarrow (\mathbb{T}^{\ast})^2$ be the map given by $$ (z_1,...,z_s,w_1,...,w_r) \mapsto Arg \left (
\frac{ \prod_{i=1}^r w^{a_i}}{ \prod_{i=1}^s z_i} , \frac{1}{\prod_{i=1}^r w_i } \right ) $$ We have:

\bigskip

\hspace{-0.6cm} \bf Proposition 3.1.2: \rm $$ lim_{ t \rightarrow -\infty} \left ( \Theta(Crit(X;f_t)) \right ) = \left \{ \left ( \frac{l \sum_{i=1}^r a_i}{(s+1)(r+1)} + \frac{k}{s+1} ,\frac {l}{r+1} \right ) \right \}_{k=0,l=0}^{s,r} \subset \mathbb{T}^2 $$

\bigskip

\hspace{-0.6cm} \bf Proof: \rm Set $A_i = lim_{t \rightarrow - \infty} A(V_i^t)$ and $B_j=lim_{t \rightarrow -\infty} A(W_j^t)$ for $1 \leq i \leq s $ and $1 \leq j \leq r$. If $(\theta , \delta) \in \bigcap_{j=1}^r B_j$ then $\delta : = \delta_1= ... = \delta_r$ and $ (r+1) \delta =0 $ in $\mathbb{T}$. If $$( \theta, \delta) \in \left (\bigcap_{i=1}^s A_i \right ) \cap \left ( \bigcap_{j=1}^r B_j \right )$$ then $ \theta=\theta_1 =... = \theta_s$ and $ \delta= \frac{l}{r+1} $ for some $0 \leq l \leq r$. As $(s+1) \theta - \sum_{j=1}^r a_j \delta$ we get $ \theta = \sum_{j=1}^r \frac{a_j l}{(s+1)(r+1)} + \frac{k}{s+1} $ for $1 \leq k \leq s$. As there are exactly $(r+1)(s+1)$ such elements $(\theta, \delta)$ we get $lim_{t \rightarrow -\infty} A(Crit(X ; f_t)) = ( \bigcap_{i=1}^s A_i ) \cap ( \bigcap_{j=1}^r B_j )$. $\square$

\bigskip

\hspace{-0.6cm} Let $i : (\mathbb{C}^{\ast})^{r+s} \rightarrow (\mathbb{C}^{\ast})^{r+s+2}$ be the map given by $$(z_1,...,z_s,w_1,...,w_r) \mapsto \left ( z_1,...,z_s,w_1,...,w_r, \frac{\prod_{i=1}^r w_i^{a_i}}{\prod_{i=1}^s z_i} , \frac{1}{\prod_{i=1}^r w_i} \right ) $$ We have:

\bigskip

\hspace{-0.6cm} \bf Corollary 3.1.4 \rm The map $E: Crit(X_a; f_t) \rightarrow Pic(X)$ given by $E(z) = [D(z)]_{\mathbb{Z}}$ is an exceptional map when $ t \rightarrow - \infty$.

\bigskip

\hspace{-0.6cm} \bf Proof: \rm Set $\Theta_{k,l}=\left ( \theta_{k,l}(a), ..., \theta_{k,l}(a), \rho_{l},...,\rho_{l}, \theta_{k,l}(a), \rho_{l} \right ) \in \mathbb{T}^{r+s+2}$ where $ \rho_l = \frac{l}{r}$ and $\theta_{k,l}(a) = \sum_{j=1}^r \frac{a_j l}{(s+1)(r+1)} + \frac{k}{s+1} $. We have: $$ [D(\Theta_{k,l})]_{\mathbb{Z}}= \left [ \sum_{i=0}^s \theta_{k,l}(a) \cdot V_X(v_i) + \sum_{j=0}^r \rho_l \cdot V_X(e_i) \right ]_{\mathbb{Z}}
= $$ $$= \left [ \left ( \frac{l \sum a_i}{r+1} + k \right ) \cdot \pi^{\ast} H + \sum_{i=0}^r \frac{ l}{r+1} \left ( \xi - a_i \cdot \pi^{\ast}H \right ) \right ]= k \cdot \pi^{\ast} H + l \cdot \xi $$ $\square$

\bigskip

\hspace{-0.6cm} Note that in the definition of the exceptional map we considered the limit $t \rightarrow - \infty$ . It is interesting to ask whether the limit $t \rightarrow \infty$ can also be interpreted in terms of the exceptional map $E$. Denote by $\Theta_{\pm}(X) :=lim_{t \rightarrow \pm} \left ( \Theta(Crit(X ; f_t)) \right )$. Consider the following example:

\bigskip

\hspace{-0.6cm} \bf Example 3.1.5 \rm (The Hirzebruch surface): Let $X = \mathbb{P}( \mathcal{O} \oplus \mathcal{O}(1))$ be the Hirzebruch surface. Recall that $$ X = \left \{ ([z_0:z_1:z_2],[\lambda_0: \lambda_1] ) \vert \lambda_0 z_0 + \lambda_1 z_1 =0 \right \} \subset \mathbb{P}^2 \times \mathbb{P}^1$$ Denote by $ p : X \rightarrow \mathbb{P}^2$ and $\pi : X \rightarrow \mathbb{P}^1$ the projection to the first and second factor, respectively. Note that $p$ expresses $X$ as the blow up of $\mathbb{P}^2$ at the point $[0:0:1] \in \mathbb{P}^2$ and $\pi$ is the fibration map. The group $Pic(X)$ is described, in turn, in the following two ways $$ Pic(X) \simeq p^{\ast} H_{\mathbb{P}^2} \cdot \mathbb{Z} \oplus E \cdot \mathbb{Z} \simeq \pi^{\ast} H_{\mathbb{P}^1} \cdot \mathbb{Z} \oplus \xi \cdot \mathbb{Z}$$ Where $E$ is class of the the line bundle whose first Chern class $c_1(E) \in H^2(X ; \mathbb{Z})$ is the Poincare dual of the exceptional divisor and $\xi$ is the class of the tautological bundle of $\pi$. The exceptional collection is expressed in these bases by $$ \mathcal{E}_X = \left \{ 0, p^{\ast} H_{ \mathbb{P}^2}-E , 2 p^{\ast} H_{\mathbb{P}^2} - E , p^{\ast} H_{\mathbb{P}^2} \right \} = \left \{ 0 , \pi^{\ast} H_{\mathbb{P}^1} , \pi^{\ast} H_{\mathbb{P}^1}+ \xi, \xi \right \} $$ Let us note that we have $p_{\ast} \left \{ 0, p^{\ast} H_{ \mathbb{P}^2} , 2 p^{\ast} H_{\mathbb{P}^2} - E \right \}= \left \{ 0, H_{\mathbb{P}^2}, 2H_{\mathbb{P}^2} \right \} = \mathcal{E}_{\mathbb{P}^2}$, while we think of the additional element $p^{\ast} H- E$ as "added by the blow up". 

\hspace{-0.6cm} On the other hand, direct computation gives $\Theta_+(X) = \mu(3) \cup \mu(1)$ where $\mu(n)=\left \{ e^{\frac{2 \pi k i}{n}} \vert k=0,...,n-1 \right \} \subset \mathbb{T} $ is the set of $n$-roots of unity for $n \in \mathbb{N}$. (see illustration in Remark 3.4). For $(k,l) \in \mathbb{Z}/2 \mathbb{Z} \oplus \mathbb{Z}/ 2 \mathbb{Z}$ let $\gamma_{kl}(t)= (z_{kl}(t),w_{kl}(t)) \in (\mathbb{C}^{\ast})^{s+r}$ be
the smooth curve defined by the condition $(z_{kl}(t) , w_{kl}(t)) \in Crit(X ; f_t)$ for $ t \in \mathbb{R}$ and $E(z_{kl}(0),w_{kl}(0)) = E_{kl}$ . Define the map $I^+: Crit(X) \rightarrow \Theta_+(X)$ by $$I^+(z_{kl},w_{kl}) :=lim_{t \rightarrow \infty}( \Theta (z_{kl}(t) ,w_{kl}(t)))$$ By direct computation $$ \begin{array}{cccccccc} I^+((z_{00},w_{00})) = \rho_3^0 & ; & I^+((z_{01},w_{01})) = \rho_3^1 & ; & I^+((z_{11},w_{11})) = \rho_3^2 & ; & I^+((z_{10},w_{10})) = 1 \end{array} $$ where $ \rho = e^{\frac{ 2 \pi i }{3} } \in \mu(3)$. Simlarly, define the map $I^-: Crit(X) \rightarrow \Theta_-(X)$,
on the other hand, taking $t \rightarrow - \infty$ in the limit. Note that this is the way we defined the exceptional map $E$ in the first place. We thus view the map $I: \Theta_-(X) \rightarrow \Theta_+(X)$ given by $I = I^+ \circ (I^-)^{-1}$ as a "geometric interpolation" between the bundle description of $\mathcal{E}_X$ and the blow up description of $\mathcal{E}_X$.

\bigskip

\subsection{Exceptional map for the class (b)} The Landau-Ginzburg potential for $X_{n,r}$ is given by $$ f(z,w) = \sum_{i=1}^{n-r} z_i + \sum_{i=1}^r w_i + \frac{1}{\prod_{i=1}^{n-r} z_i} + \frac{1}{\prod_{i=1}^r w_i} + \frac{1}{\prod_{i=1}^{n-r} z_i \cdot \prod_{i=1}^r w_i } \in L (\Delta_{n,r}^{\circ})$$ Instead, for $t \in \mathbb{R}$  we consider the $1$-parametric family $$f_{t} (z,w) = \sum_{i=1}^{n-r} e^{-t} z_i + \sum_{i=1}^r e^{-t} w_i + \frac{1}{\prod_{i=1}^{n-r} z_i} + \frac{1}{\prod_{i=1}^r w_i} + \frac{1}{\prod_{i=1}^{n-r} z_i \cdot \prod_{i=1}^r w_i } \in L (\Delta_{n,r}^{\circ})$$ The corresponding system of equations is $$ \begin{array}{ccc} e^{-t} z_i - \frac{1}{\prod_{i=1}^{n-r} z_i}- \frac{1}{\prod_{i=1}^{n-r} z_i \cdot \prod_{i=1}^r w_i } =0 & ; &
e^{-t} w_j - \frac{1}{\prod_{i=1}^{r} w_i}- \frac{1}{\prod_{i=1}^{n-r} z_i \cdot \prod_{i=1}^r w_i } =0 \end{array}$$ for $1 \leq i \leq n-r$ and $1 \leq j \leq r$. One can see that $ z_1 =...= z_{n-r}$ and $ w_1 = ... = w_r $, and we set $z=z_i$ and $w= w_j$.
In particular, the equations turn to be $$ \left \{ \begin{array}{c} e^{-t} z - \frac{1}{z^{n-r}}- \frac{1}{z^{n-r} w^r } =0 \\
e^{-t} w - \frac{1}{w^r}- \frac{1}{z^{n-r} w^r } =0 \end{array} \right.  $$ Let $\Theta: (\mathbb{C}^{\ast})^n \rightarrow \mathbb{T}^2$ given by $(z_1,...,z_{n-r},w_1,...,w_r) \mapsto Arg(z_1,w_1)$. We have:

\bigskip

\hspace{-0.6cm} \bf Proposition 3.2.1: \rm The elements of $lim_{t \rightarrow \infty} \Theta(Crit(X_{n,r}; f_t)) \subset \mathbb{T}^{2}$ are of the following two classes $$ \Theta_{k,l}= \left (k \cdot \rho_{n-r+1}, l \cdot \rho_{r+1} \right ) \hspace{0.5cm} \textrm{ for } \hspace{0.25cm} 0 \leq k \leq n-r, 0 \leq l \leq r$$ and $$ \Theta'_{m,n} = (\rho_{2(n-r)} \cdot m \cdot \rho_{n-r} , \rho_{2r} \cdot n \cdot \rho_{r} ) \hspace{0.5cm} \textrm{ for } \hspace{0.25cm} 0 \leq m \leq n-r-1, 0 \leq n \leq r-1$$

\hspace{-0.6cm} \bf Proof: \rm Geometrically, the elements of $Crit(X_{n,r},f_t) \subset (\mathbb{C}^{\ast})^2$ are given as the intersection of the two curves $$\begin{array}{ccc} C_1 := \left \{ e^{-t} z - \frac{1}{z^{n-r}} - \frac{1}{z^{n-r} \cdot w^r} =0 \right \} \subset (\mathbb{C}^{\ast})^2 & ; & C_2 := \left \{ e^{-t} w - \frac{1}{w^{r}} - \frac{1}{z^{n-r} \cdot w^r} =0 \right \} \subset (\mathbb{C}^{\ast})^2 \end{array}$$ Let $Log_r \vert \cdot \vert : (\mathbb{C}^2) \rightarrow \mathbb{R}^2$ and $Arg: (\mathbb{C}^{\ast})^2 \rightarrow \mathbb{T}^2$ be the maps given by $$ \begin{array}{ccc} (r_1 e^{2 \pi i \theta_1} , r_2 e^{2 \pi i \theta_2}) \mapsto\left ( \frac{log(r_1)}{log(r)} , \frac{log (r_2)}{log(r)} \right ) & ; & (r_1 e^{2 \pi i \theta_1} , r_2 e^{2 \pi i \theta_2}) \mapsto (\theta_1 ,\theta_2 ) \end{array}$$ respectively. Let $T(C_i) := Log \vert C_1 \vert \subset \mathbb{R}^2 $ and $A(C_i) := Arg(C_i)$ be the amobea and co-amobea of $C_i$ for $i=1,2$. The corresponding tropical curves $\widetilde{T}(C_i) := lim_{r \rightarrow \infty} Log_r \vert C_i \vert$ are given as follows: $$ T(C_1) = \left \{ \begin{array}{c} x_1 =\frac{t}{n-r+1} \textrm{ , } \\ x_2 \leq 0 \end{array} \right \} \cup \left \{\begin{array}{c} x_2 = 0 \textrm{ , } \\ x_1 \leq \frac{t}{n-r+1} \end{array} \right \} \cup \left \{ \begin{array}{c} (n-r+1) x_1+ r x_2 = t \textrm{ , } \\ x_2 \leq 0 \end{array}
\right \} $$
$$ T(C_2) = \left \{ \begin{array}{c} x_2 =\frac{t}{r+1} \textrm{ , } \\ x_2 \leq 0 \end{array} \right \} \cup \left \{\begin{array}{c} x_1 = 0 \textrm{ , } \\ x_2 \leq \frac{t}{r+1} \end{array} \right \} \cup \left \{ \begin{array}{c} (n-r) x_1+ (r+1) x_2 = t \textrm{ , } \\ x_1 \leq 0 \end{array}
\right \} $$ These two tropical curves, which are the spines of the corresponding amoeba, intersect at two points $$ \widetilde{T}(C_1 ) \cap \widetilde{T}(C_2) = \left \{ ( 0,0) , \left ( \frac{t}{n-r+1}, \frac{t}{r+1} \right ) \right \} \subset \mathbb{R}^2$$ The corresponding co-tropical curves $\widetilde{A}(C_i) := \partial A(C_i)$ are given by: $$ A(C_1) = \left \{ \theta_1 =\frac{k}{n-r+1} \right \}_{k=0}^{n-r} \cup \left \{ r \theta_2 =\frac{1}{2}+k \right \}_{k=0}^{r-1} \cup \left \{ (n-r+1) \theta_1+ r \theta_2 = k \right \}_{k=0}^{r-1} $$ $$ A(C_2) = \left \{ \theta_2 =\frac{k}{r+1} \right \}_{k=0}^{r+1} \cup \left \{ (n-r) \theta_1 =\frac{1}{2}+k \right \}_{k=0}^{n-r-1} \cup \left \{ (n-r) \theta_1+ (r+1) \theta_2 = k \right \}_{k=0}^{n-r-1} $$ Each component of the co-tropical curve corresponds to one of the tentecals of the tropical curve. Thus, for the intersection point of the tropical curves $$\left ( \frac{t}{n-r+1}, \frac{t}{r+1} \right ) \in \widetilde{T}(C_1) \cap \widetilde{T}(C_2)$$ the corresponding intersections for the co-tropical curves are $$ \left \{ (\theta_1, \theta_2) \vert \theta_1 = \frac{k_1}{n-r+1} \textrm{ and } \theta_2 = \frac{k_2}{r+1} \right \}_{k_1=0,k_2=0}^{n-r,r} \subset \mathbb{T}^2 $$ For the intersection point of the tropical curves $ \left ( 0,0 \right ) \in \widetilde{T}(C_1) \cap \widetilde{T}(C_2)$ the corresponding intersections for the co-tropical curves are $$ \left \{ (\theta_1, \theta_2) \vert \theta_1 = \frac{1}{2(n-r)} +\frac{k_1}{n-r} \textrm{ and } \theta_2 = \frac{1}{2r} +\frac{k_2}{r} \right \}_{k_1=0,k_2=0}^{n-r-1,r-1} \subset \mathbb{T}^2 $$ $\square$

\bigskip

\hspace{-0.6cm} We have:

\bigskip

\hspace{-0.6cm} \bf Corollary 3.2.2: \rm The map $E: Crit(X_{n,r} ; f_t) \rightarrow Pic(X)$ given by $E(z) =[D(z)]_{\mathbb{Z}}$ is an exceptional map when $ t \rightarrow \infty$.

\bigskip

\hspace{-0.6cm} \bf Proof: \rm We have $$ [D(\Theta_{k,l})]_{\mathbb{Z}} = \left [ \frac{ (n-r) k}{n-r+1} \cdot U+ \frac{r \cdot l}{r+1} \cdot V + \right. $$ $$ \left.+\frac{k}{n-r+1} \left (U-E \right ) + \frac{l}{r+1} \left ( V-E \right ) + \left ( \frac{k}{n-r+1} + \frac{l}{r+1} \right ) E \right ]_{\mathbb{Z}} = k \cdot U + l \cdot V$$ and
$$ [D (\Theta'_{m,n})]_{\mathbb{Z}} = \left [\left ( \frac{1}{2} +m \right ) \cdot U+ \left ( \frac{1}{2} + n \right ) \cdot V + \frac{1}{2} \left (U-E \right ) + \frac{1}{2} \left ( V-E \right ) \right ]_{\mathbb{Z}} = (1+m) \cdot U +(1+n) \cdot V - E $$ $\square$

\bigskip

\subsection{Exceptional map for the class (c)} The Landau-Ginzburg potential for $X_{n,b}$ is given by $$ f(z,w) = \sum_{i=1}^{n-1} z_i +e^{-t} w + \frac{1}{\prod_{i=1}^{n-1} z_i \cdot w^b} + \frac{1}{\prod_{i=1}^{n-1} z_i \cdot w^{b+1}} + \frac{1}{ w } \in L (\Delta_{n,b}^{\circ})$$ Instead, for $t \in \mathbb{R}$ we consider the $1$-parametric family$$ f_t(z,w) = \sum_{i=1}^{n-1} e^{-t} z_i + e^{-t} w + \frac{1}{\prod_{i=1}^{n-1} z_i \cdot w^b} -\frac{i}{\prod_{i=1}^{n-1} z_i \cdot w^{b+1}} + \frac{1}{ w } \in L (\Delta_{n,b}^{\circ})$$ The corresponding system of equations is $$ \begin{array}{ccc} e^{-t} z_i - \frac{1}{\prod_{i=1}^{n-1} z_i \cdot w^b}+\frac{i}{\prod_{i=1}^{n-1} z_i \cdot w^{b+1} } =0 & ; &
e^{-t} w - \frac{b}{\prod_{i=1}^{n-1} z_i \cdot w^b}+ \frac{(b+1)i}{\prod_{i=1}^{n-1} z_i \cdot w^{b+1} } - \frac{1}{w} =0 \end{array}$$ for $1 \leq i \leq n-1$. One can see that $ z_1 =...= z_{n-1}$ and we set $z=z_i$.
In particular, the equations turn to be $$ \left \{ \begin{array}{c} e^{-t} z - \frac{1}{z^{n-1} \cdot w^b}+\frac{i}{z^{n-1} w^{b+1} } =0 \\ e^{-t} w - \frac{b}{z^{n-1} \cdot w^b}+ \frac{(b+1)i}{z^{n-1} w^{b+1} }- \frac{1}{w} =0 \end{array} \right. $$ Let $\Theta: (\mathbb{C}^{\ast})^n \rightarrow \mathbb{T}^2$ given by $(z_1,...,z_{n-1},w) \mapsto Arg(z_1,w)$. We have:

\bigskip

\hspace{-0.6cm} \bf Proposition 3.3.1: \rm The elements of $lim_{t \rightarrow \infty} \Theta(Crit(X_{n,b}; f_t)) \subset \mathbb{T}^{2}$ are of the following two classes $$ \Theta_{k,l} = \left (( b \cdot l \cdot \rho_{2n}) \cdot (k \cdot \rho_n), \frac{l}{2} \right ) \hspace{0.5cm} \textrm{ for } \hspace{0.25cm} 0 \leq k \leq n-1, 0 \leq l \leq 1$$ and $$ \Theta'_{m} = \left ( (3(b+1) \cdot \rho_{4(n-1)} ) \cdot (m \cdot \rho_{n-1}) , \frac{1}{4} \right ) \hspace{0.5cm} \textrm{ for } \hspace{0.25cm} 0 \leq m \leq n-2$$
\bigskip

\hspace{-0.6cm} \bf Proof: \rm As $ t \rightarrow \infty $ the solutions converge to the union of the solutions of the two systems $$ \left \{ \begin{array}{c} - \frac{1}{z^{n-1} \cdot w^b}+\frac{i}{z^{n-1} w^{b+1} } =0 \\ - \frac{b}{z^{n-1} \cdot w^b}+ \frac{(b+1)i}{z^{n-1} w^{b+1} }- \frac{1}{w} =0 \end{array} \right. \hspace{0.25cm} \textrm{ and } \hspace{0.25cm} \left \{ \begin{array}{c} e^{-t} z - \frac{1}{z^{n-1} \cdot w^b} =0 \\ e^{-t} w - \frac{b}{z^{n-1} \cdot w^b}- \frac{1}{w} =0 \end{array} \right. $$ The solutions of the first system are given by direct computation by $$ \left \{ (\rho_{4(n-1)}^{3(b+1)} \cdot \rho_{n-1}^m, i) \vert m=0,...,n-2 \right \} $$ For the second equation we have $$ e^{-t} w -be^{-t} z - \frac{1}{w} = 0 $$ solving this quadratic equation gives $$ w_{\pm} = \frac{ be^{-t} z \pm \sqrt{b z e^{-2t} + 4 e^{-t} } }{ 2e^{-t}} \rightarrow \pm e^{\frac{t}{2}} $$ Substituting in the first equation gives $ z^n w^b e^{-t} = 1 $. For $w_+$ we have $$ z= e^{\frac{(2-b)t}{2n}} \rho_n^k \hspace{0.5cm} \textrm{ for } \hspace{0.25cm} k=0,...,n-1$$ For $w_-$ we have $$ z= e^{\frac{(2-b)t}{2n}}\rho^b_{2n} \cdot \rho_n^k \hspace{0.5cm} \textrm{ for } \hspace{0.25cm} k=0,...,n-1$$ which gives the required result. $\square$

\bigskip

\hspace{-0.6cm} We have:

\bigskip

\hspace{-0.6cm} \bf Corollary 3.3.2: \rm The map $E: Crit(X_{n,b} ; f_t) \rightarrow Pic(X)$ given by $E(z) =[D(z)]_{\mathbb{Z}}$ is an exceptional map when $ t \rightarrow \infty$.

\bigskip

\hspace{-0.6cm} \bf Proof: \rm For $0 \leq k \leq n-1$ we have $$ [D(\Theta_{k,0})]_{\mathbb{Z}}=\left [ \frac{(n-1)k}{n} \cdot V + \frac{k}{n} \cdot (V-Y) + \left ( \frac{k}{n} + \frac{3}{4} \right ) Y \right]_{\mathbb{Z}} = k \cdot V $$ and $$ [D(\Theta_{k,1})]_{\mathbb{Z}} =\left [ \left ( \frac{(n-1)b}{2n}+\frac{(n-1)k}{n} \right ) \cdot V + \frac{1}{2}(Y+T+bV) + \left (\frac{b}{2n}+\frac{k}{n} \right) \cdot (V-Y) \right.$$ $$ \left. +\left ( \frac{b}{2n} + \frac{k}{n} + \frac{1}{2}+\frac{3}{4} \right ) Y + \frac{1}{2}T \right ]_{\mathbb{Z}} = k \cdot V + (Y+T+bV)$$ and $$ [D( \Theta'_{m,n})]_{\mathbb{Z}} = \left [
\left ( \frac{3b+1}{4}+m \right ) \cdot V + \frac{1}{4} (Y+T+bV) + \frac{1}{4} (V-Y) \right ]_{\mathbb{Z}} = (m+1) \cdot V + T $$ for $0 \leq m \leq n-2$. $\square$

\bigskip

\hspace{-0.6cm} Let us conclude this section with the following remark:

\bigskip

\hspace{-0.6cm} \bf Remark 3.3.3 \rm (Exceptional maps and the Frobenius mapping): Let $\widetilde{F} : \mathbb{T}^{r} \rightarrow Pic(X)$, where $r = rk(Div_T(X))$, be the Frobenius map described in Remark 2.6. Consider the map $i : (\mathbb{C}^{\ast})^n \rightarrow (\mathbb{C}^{\ast})^r$ given by $z \mapsto (z^{n_1},...,z^{n_r})$ where $\Delta^{\circ}(0)= \left \{ n_1,...,n_r \right \}$. Note that by definition of the Frobenius map, we can equivalently write the exceptional maps as $E(z) = \widetilde{F}(Arg(i(z)))$ for $z \in Crit(X ; f_t)$ for $0<<t$ big enough. In particular, we geometrically interpret the torus $\mathbb{T}^r$ (which in the Frobenius setting was defined as the limit of $Div_T(X)/l \cdot Div_T(X)$ as $l \rightarrow \infty$) as the argument torus of the maximal algebraic tori of $\mathbb{P}(L(\Delta^{\circ})^{\vee})$.

\section{Monodromies and the Endomorphism Ring}
\label{s:LGM}

\hspace{-0.6cm} Given a full strongly exceptional collection $\mathcal{E} = \left \{ E_i \right \}_{i=1}^N \subset Pic(X)$ one is interested in the structure of its endomorphism algebra $$A_{\mathcal{E}}= End \left ( \bigoplus_{i=1}^N E_i \right )= \bigoplus_{i,j=0}^N Hom(E_i,E_j)= \bigoplus_{i,j=0}^N H^0(X ; E_j \otimes E_i^{-1})$$ Our aim in this section is to illustrate how this algebra is naturally reflected in the monodromy group action of the Landau-Ginzburg system and the Frobenius stratification, for the classes of toric Fano manifolds (a)-(c), considered above.

\hspace{-0.6cm} Recall that a \emph{quiver with relations} $\widetilde{Q}=(Q,R)$ is a directed graph $Q$ with a two sided ideal $R$ in the path algebra $\mathbb{C}Q$ of $Q$, see \cite{DW}. In particular, a quiver with relations $\widetilde{Q}$ determines the
associative algebra $A_{\widetilde{Q}}=\mathbb{C}Q/R$, called the path algebra of $\widetilde{Q}$. In general, a collection of elements $\mathcal{C} \subset \mathcal{D}^b(X)$ and a basis $ B \subset A_{\mathcal{C}}:= End \left ( \bigoplus_{E \in \mathcal{C}} E \right )$
determine a quiver with relations
$\widetilde{Q}(\mathcal{C},B)$ whose vertex set is $\mathcal{C}$ such that $A_{\mathcal{C}} \simeq A_{Q(\mathcal{C},B)}$, see \cite{K}.

\hspace{-0.6cm} Let $Gr(X)$ be the Grassmaniann of subspaces $W \subset \mathbb{P}(L(\Delta^{\circ})^{\vee})$ of $codim(W)=n$. Denote by $Crit(X ; W) := X^{\circ} \cap W$ for $W \in Gr(X)$. Note the natural map $ W : L(\Delta^{\circ}) \rightarrow Gr(X)$ given by $$ f \mapsto \bigcap_{i=1}^n \left \{ z_i \frac{\partial}{\partial z_i} f =0 \right \} $$ In particular $Crit(X ; W(f)) = Crit(X ; f)$. Let $R_X \subset Gr(X)$ be the hypersurface of all $W \in Gr(X)$ such that $Crit(X;W)$ is non-reduced. Whenever $Crit(X; W(f_u))$ is reduced, one obtains, via standard analytic continuation, a monodromy map of the following form $$ M : \pi_1(Gr(X)
\setminus R_X, W(f_u)) \rightarrow Aut(Crit(X; W(f_u)))$$ On the other hand, as the line bundles in the exceptional collections $ \mathcal{E}$ are all $T$-equivariant the $Hom$-spaces between the elements of the collection admit a decomposition $$ Hom(E_i, E_j) = \bigoplus_{D \in Div_T(X)} Hom_D(E_i,E_j)$$ For each of the classes (a)-(c) we define a natural map $\Gamma : Div_T(X) \rightarrow \pi_1(Gr(X) \setminus R_X , f_u)$ with the property $$ Hom_D(E(z_i),E(z_j)) \neq 0 \Rightarrow M(\Gamma_D)(z_i)=z_j $$ We thus define $$ Hom_{mon} (z_i , z_j) := \bigoplus_{\left \{ D \vert M(\Gamma_D)(z_i)=z_j \right \}} D \cdot \mathbb{Z} $$ Our aim in this section is to show how exceptional maps $E : Crit(X ; f_u ) \rightarrow Pic(X)$ of section 3 satisfy $$Hom(E(z), E(w) ) \subset Hom_{mon} (z,w) \hspace{0.5cm} \textrm{ for any } \hspace{0.25cm} z,w \in Crit(X ; f_u)$$

\bigskip

\hspace{-0.6cm} \bf Remark 4.1: \rm We refer to this as the $M$-aligned property to note that the exceptional maps in section 3 are defined such that the algebraic $Hom$-functor is aligned with the geometric monodromies of the Landau-Ginzburg system. Note that this property is non-trivial as, in particular, it implies that even for an individual $T$-divisor $D \in Div_T(X)$, the corresponding monodromy $M(\Gamma_D)$ is aligned for all elements of the collection at once.

\subsection{M-algined property for projective bundles} Let $X= \mathbb{P} \left (\mathcal{O}_{\mathbb{P}^s} \oplus \bigoplus_{i=1}^r \mathcal{O}_{\mathbb{P}^s}(a_i) \right ) $ be a projective Fano bundle. Note that, in our case $$Div_T(X) = \left ( \bigoplus_{i=1}^s \mathbb{Z} \cdot V_X(v_i) \right )
\bigoplus \left ( \bigoplus_{i=0}^r \mathbb{Z} \cdot V_X(e_i) \right )$$ Let $L_{kl}=k \cdot \pi^{\ast} H + l \cdot \xi \in Pic(X)$ be any element. Then $$ H^0 (X ; L_{kl}) \simeq \left \{ \sum_{i=0}^s n_i V_X(v_i) + \sum_{i=0}^r m_i V_X(e_i) \bigg \vert \vert m \vert = l \textrm{ and } \vert n \vert = k+ \sum_{i=0}^r m_i a_i \right \} \subset Div^+_T(X)$$

\hspace{-0.6cm} In particular, the algebra $A_{\mathcal{E}}$ comes with the basis $\left \{ V(v_0),...,V(v_s),V(e_0),...,V(e_r) \right \}$. We denote the resulting quiver by $Q_s(a_0,...,a_r)$. For example, the quiver $Q_3(0,1,2)$ for $X=\mathbb{P}\left ( \mathcal{O}_{\mathbb{P}^3} \oplus \mathcal{O}_{\mathbb{P}^3}(1) \oplus \mathcal{O}_{\mathbb{P}^3}(2) \right )$ is the following

$$
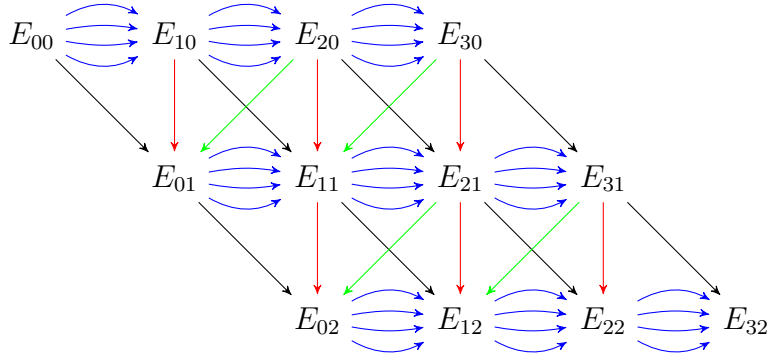
\begin{tikzpicture}[->,>=stealth',shorten >=1pt,auto,node distance=1.9cm,main node/.style={font=\sffamily\bfseries \small}]
\node (1) {$E_{00}$};
\node (4) [right of=1] {$E_{10}$};
\node (7) [right of=4] {$E_{20}$};
\node (10) [right of=7] {$E_{30}$};

\node (2) [below of=4] {$E_{01}$};
\node (5) [right of=2] {$E_{11}$};
\node (8) [right of=5] {$E_{21}$};
\node (11) [right of=8] {$E_{31}$};

\node (3) [below of=5] {$E_{02}$};
\node (6) [right of=3] {$E_{12}$};
\node (9) [right of=6] {$E_{22}$};
\node (12) [right of=9] {$E_{32}$};

\path[every node/.style={font=\sffamily\small}]
(1) edge [bend right=30,blue] node {} (4)
edge [bend right =10,blue] node {} (4)
edge [bend left=30, blue] node {} (4)
edge [bend left = 10,blue] node {} (4)
edge node {} (2)

(2) edge [bend right=30,blue] node {} (5)
edge [bend right = 10,blue] node {} (5)
edge [bend left=30, blue] node {} (5)
edge [bend left = 10, blue] node{} (5)
edge node {}(3)
(3) edge [bend right=30, blue] node {} (6)
edge [bend right = 10,blue] node {} (6)
edge [bend left=30, blue] node {} (6)
edge [bend left = 10,blue] node{} (6)
(4) edge [bend right=30, blue] node {} (7)
edge [bend right = 10,blue] node {} (7)
edge [bend left=30, blue] node {} (7)
edge [bend left = 10, blue] node{} (7)
edge node {}(5)
edge [red] node {}(2)

(5) edge [bend right=30, blue] node {} (8)
edge [bend right = 10,blue] node {} (8)
edge [bend left=30, blue] node {} (8)
edge [bend left =10, blue] node{} (8)
edge node {}(6)
edge [red] node {}(3)
(6) edge [bend right, blue] node {} (9)
edge [bend right =10,blue] node {} (9)
edge [bend left=30, blue] node {} (9)
edge [bend left=10,blue] node{} (9)
(7) edge node {}(8)
edge [red] node {}(5)
edge [green] node {}(2)
edge [bend right=30, blue] node {} (10)
edge [bend right=10,blue] node {} (10)
edge [bend left=30, blue] node {} (10)
edge [bend left=10,blue] node{} (10)
(8) edge node {}(9)
edge [red] node {}(6)
edge [green] node {}(3)
edge [bend right=30, blue] node {} (11)
edge [bend left=30,blue] node {} (11)
edge [bend left=10, blue] node {} (11)
edge [bend right=10,blue] node{} (11)
(9) edge [bend right=30, blue] node {} (12)
edge [bend right=10,blue] node {} (12)
edge [bend left=30, blue] node {} (12)
edge [bend left =10,blue] node{} (12)
(10) edge node{} (11)
edge [red] node{} (8)
edge [green] node{} (5)
(11) edge node{} (12)
edge [red] node{} (9)
edge [green] node{} (6);
\end{tikzpicture}
$$ For a divisor $D= \sum_{i=0}^s n_i V_X(v_i) + \sum_{i=0}^r m_i V_X(e_i) \in Div_T(X) $ and $u \in \mathbb{C}$ consider the loop $$ \gamma^u_{D}(\theta) := \sum_{i=1}^s e^{2 \pi i n_i \theta } z_i + \sum_{i=1}^{r} e^{2 \pi i m_i \theta} w_i + e^{u} \cdot e^{2 \pi i n_0 \theta }
\frac{ \prod_{i=1}^r w_i^{a_i}}{ \prod_{i=1}^s z_s} + \frac{e^{2 \pi i m_0 \theta} }{\prod_{i=1}^r w_i } $$ For $ \theta \in [0,1)$. Define $ \Gamma_D : = lim_{t \rightarrow -\infty} [\gamma_D^t] \in \pi_1(L(\Delta^{\circ}) \setminus R_X , f_X) $ and set $\widetilde{M}_D := M(\Gamma_D) \in Aut(Crit(X;f_t))$. Express the solution scheme as $$Crit(X;f_t) = \left \{ (z_{kl},w_{kl} ) \right \}_{k=0,l=0}^{s,r} \simeq \mathbb{Z} / (r+1) \mathbb{Z} \oplus \mathbb{Z} / (s+1) \mathbb{Z}$$ where $E((z_{kl},w_{kl}))=E_{kl}$. We have:

\bigskip

\hspace{-0.6cm} \bf Theorem 4.2 \rm (M-aligned property (a)): For $(k,l) \in \mathbb{Z}/(s+1) \mathbb{Z} \oplus \mathbb{Z} / (r+1) \mathbb{Z} \simeq Crit(X)$ the monodromy action satisfies:
\bigskip

(a) $\widetilde{M}_{V(v_j)}(k,l)=(k+1,l)$ for $j=0,...,s$.

\bigskip

(b) $\widetilde{M}_{V(e_j)}(k,l)=(k-a_j,l+1)$ for $j=0,...,r$.

\bigskip

\hspace{-0.6cm} \bf Proof : \rm For a divisor $D \in Div_T(X)$ and $ \theta \in [0,1)$ Set $$ \begin{array}{ccc} V^{u, \theta}_{D,i} := \left \{ e^{2 \pi i n_i \theta} z_i - e^u e^{ 2 \pi i n_0 \theta} \frac{ \prod_{i=1}^r w_i^{a_i}}{\prod_{i=1}^s z_i }=0 \right \} & ; & W^{u , \theta}_{D,j} := \left \{ e^{2 \pi i m_j}w_j +a_i e^u e^{2 \pi i n_0 \theta} \frac{ \prod_{i=1}^r w_i^{a_i}}{ \prod_{i=1}^s z_i} - \frac{e^{ 2 \pi i m_0}}{\prod_{i=1}^r w_i } = 0 \right \} \end{array}$$ where $ 1 \leq i \leq s$, $1 \leq j \leq r$ and $u \in \mathbb{C}$. Let $( \theta_1,...,\theta_s, \delta_1,...,\delta_r)$ be coordinates on $\mathbb{T}^{s+r}$. It is clear that:

\bigskip

- $A_{D,i}^{ t, \theta}:= lim_{ t \rightarrow -\infty} A(V^{t,\theta}_{D,i}) = \left \{ \theta_i + \sum_{i=1}^s \theta_i - \sum_{j=1}^r a_j \delta_j +(n_i-n_0) \theta=0 \right \} \subset \mathbb{T}^{s+r}$

\bigskip

- $B_{D,j}^{t , \theta}:= lim_{t \rightarrow -\infty} A(W_{D,j}^{t, \theta})= \left \{ \delta_j + \sum_{j=1}^r \delta_j + (m_j-m_0) \theta =0 \right \} \subset \mathbb{T}^{s+r}$

\bigskip

\hspace{-0.6cm} For $D = V(v_0) $ we have $ (\theta,\delta) \in \bigcap_{j=1}^r B_{D,j}^{t,\theta}$ then $ \delta: = \delta_1 = ... = \delta_r$ and $ (r+1) \delta = 0$ hence $ \delta = \frac{l}{r+1}$ for some $0 \leq l \leq r$. Assume further that $ (\theta,\delta) \in (\bigcap_{i=1}^s A_{D,i}^{t, \theta}) \cap( \bigcap_{j=1}^r B_{D,j}^{t,\theta})$ then $ \widetilde{ \theta} = \theta_1 = ... = \theta_s$ and $ (s+1) \widetilde{ \theta} - \sum_{j=1}^r \frac{ a_j l }{r+1} - \theta= 0$. Hence, $\widetilde{\theta} = \frac{k}{s+1} + \frac{ l \sum_{j=1}^r a_j}{(s+1)(r+1)} + \frac{\theta}{(s+1)} $ for some $0 \leq k \leq s$.

\hspace{-0.6cm} For $D= V(v_i)$ if $ (\theta,\delta) \in (\bigcap_{i=1}^s A_{D,i}^{t, \theta}) \cap( \bigcap_{j=1}^r B_{D,j}^{t,\theta})$ then $ \widetilde{ \theta} = \theta_1 = ...= \hat{\theta_i}=... = \theta_s$ and $ \theta_i = \widetilde{ \theta} - \theta$ and again $ (s+1) \widetilde{ \theta} - \sum_{j=1}^r \frac{ a_j l }{r+1} - \theta= 0$. Hence, $\widetilde{\theta} = \frac{k}{s+1} + \frac{ l \sum_{j=1}^r a_j}{(s+1)(r+1)} + \frac{\theta}{(s+1)} $ for some $0 \leq k \leq s$

\hspace{-0.6cm} For $D=V(e_0)$ we have $ (\theta,\delta) \in \bigcap_{j=1}^r B_{D,j}^{t,\theta}$ then $ \delta: = \delta_1 = ... = \delta_r$ and $ (r+1) \delta = \theta $ hence $ \delta = \frac{l+ \theta}{r+1}$ for some $0 \leq l \leq r$. Assume further that $ (\theta,\delta) \in (\bigcap_{i=1}^s A_{D,i}^{t, \theta}) \cap( \bigcap_{j=1}^r B_{D,j}^{t,\theta})$ then $ \widetilde{ \theta} = \theta_1 = ... = \theta_s$ and $ (s+1) \widetilde{ \theta} - \sum_{j=1}^r \frac{ a_j (l+ \theta) }{r+1} = 0$. Hence, $\widetilde{\theta} = \frac{k}{s+1} + \frac{ (l+ \theta) \sum_{j=1}^r a_j}{(s+1)(r+1)} $ for some $ 0 \leq k \leq s$.

\hspace{-0.6cm} For $D=V(e_j)$ we have $ (\theta,\delta) \in \bigcap_{j=1}^r B_{D,j}^{t,\theta}$ then $ \delta: = \delta_1 = ... = \hat{ \delta_j} =...= \delta_r$ and $\delta_j= \delta - \theta$ hence $(r+1) \delta = \theta $ and $ \delta = \frac{l+ \theta}{r+1}$ for some $0 \leq l \leq r$. Assume $ (\theta,\delta) \in (\bigcap_{i=1}^s A_{D,i}^{t, \theta}) \cap( \bigcap_{j=1}^r B_{D,j}^{t,\theta})$ then $ \widetilde{ \theta} = \theta_1 = ... = \theta_s$ and $ (s+1) \widetilde{ \theta} - \sum_{j=1}^r \frac{ a_j (l+ \theta) }{r+1} + a_j \theta= 0$. Hence, $\widetilde{\theta} = \frac{k-a_j \theta}{s+1} + \frac{ (l+ \theta) \sum_{j=1}^r a_j}{(s+1)(r+1)} $ for some $ 0 \leq k \leq s$. $ \square$

\bigskip

\hspace{-0.6cm} For instance, consider the following example:

\bigskip

\hspace{-0.6cm} \bf Example 4.3 \rm (monodromies for $X=\mathbb{P}\left ( \mathcal{O}_{\mathbb{P}^3} \oplus \mathcal{O}_{\mathbb{P}^3}(1) \oplus \mathcal{O}_{\mathbb{P}^3}(2) \right )$): The following diagram outlines the corresponding monodromies
on $\mathbb{T}^2$:

$$
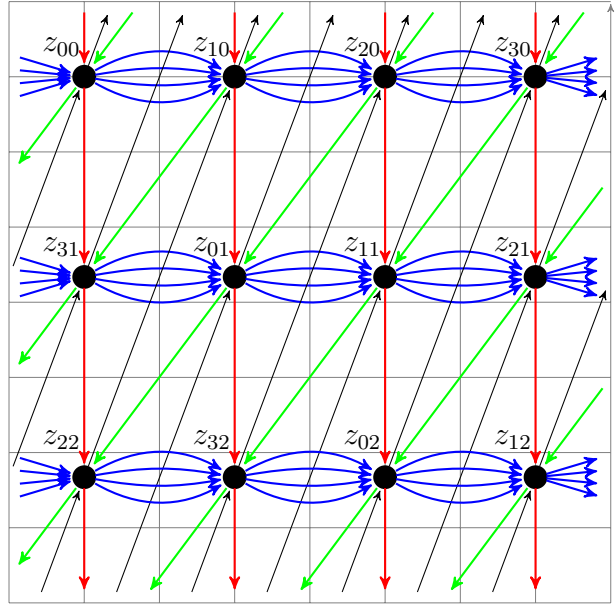
\begin{tikzpicture}[->,>=stealth',shorten >=1pt,auto,node distance=1.3cm,main node/.style={circle,draw,font=\sffamily\bfseries }, scale=1]
\draw[step=1cm,gray,very thin] (-4,-8) grid (4,0);
\node (a0) at (-1,0){};
\node (b0) at (1,0){};
\node (c0) at (3,0){};
\node (d0) at (-3,0){};

\node (a05) at (-0.249,0){};
\node (b05) at (1.75,0){};
\node (c05) at (3.75,0){};
\node (d05) at (-2.249,0){};

\node (a55) at (-0.249,-8){};
\node (b55) at (1.75,-8){};
\node (c55) at (3.75,-8){};
\node (d55) at (-2.249,-8){};

\node (a00) at (-0.624,0){};
\node (b00) at (1.375,0){};
\node (c00) at (3.375,0){};
\node (d00) at (-2.624,0){};

\node (a01) at (0.375,0){};
\node (b01) at (2.375,0){};

\node (e00) at (4,-1){};
\node (d01) at (-1.624,0){};

\node (a44) at (-2.624,-8){};
\node (b44) at (-1.624,-8){};
\node (c44) at (-0.624,-8){};
\node (d44) at (0.375,-8){};
\node (e44) at (1.375,-8){};
\node (f44) at (2.375,-8){};
\node (g44) at (-4,-6.328){};
\node (h44) at (-3.624,-8){};

\node (c5) at (4,-5){};
\node (d5) at (-4,-5){};

\node (c6) at (4,-2.33){};
\node (d6) at (-4,-2.33){};

\node (c7) at (4,-3.664){};
\node (d7) at (-4,-3.664){};

\node (d8) at (-4,-7.664) {};

\node (a1) at (-1,-1){};
\node (b1) at (1,-1){};
\node (c1) at (3,-1){};
\node (d1) at (-3,-1){};
\node (c11) at (4,-1.3){};
\node (c12) at (4,-1.1){};
\node (c13) at (4,-0.7){};
\node (c14) at (4,-0.9){};
\node (d11) at (-4,-1.3){};
\node (d12) at (-4,-1.1){};
\node (d13) at (-4,-0.7){};
\node (d14) at (-4,-0.9){};

\node (a2) at (-1,-3.664){};
\node (b2) at (1,-3.664){};
\node (c2) at (3,-3.664){};
\node (d2) at (-3,-3.664){};
\node (c21) at (4,-3.964){};
\node (c22) at (4,-3.764){};
\node (c23) at (4,-3.364){};
\node (c24) at (4,-3.564){};
\node (d21) at (-4,-3.964){};
\node (d22) at (-4,-3.764){};
\node (d23) at (-4,-3.364){};
\node (d24) at (-4,-3.564){};

\node (a3) at (-1,-6.328){};
\node (b3) at (1,-6.328){};
\node (c3) at (3,-6.328){};
\node (d3) at (-3,-6.328){};
\node (c31) at (4,-6.628){};
\node (c32) at (4,-6.428){};
\node (c33) at (4,-6.028){};
\node (c34) at (4,-6.228){};
\node (d31) at (-4,-6.628){};
\node (d32) at (-4,-6.428){};
\node (d33) at (-4,-6.028){};
\node (d34) at (-4,-6.228){};

\node (a4) at (-1,-8){};
\node (b4) at (1,-8){};
\node (c4) at (3,-8){};
\node (d4) at (-3,-8){};

\node (aa1) at (-3.3,-0.6){$z_{00}$};
\node (bb1) at (-1.3,-0.6){$z_{10}$};
\node (cc1) at (0.7,-0.6){$z_{20}$};
\node (dd1) at (2.7,-0.6){$z_{30}$};

\node (aa2) at (-1.3,-3.264){$z_{01}$};
\node (bb2) at (0.7,-3.264){$z_{11}$};
\node (cc2) at (2.7,-3.264){$z_{21}$};
\node (dd2) at (-3.3,-3.264){$z_{31}$};

\node (aa3) at (-1.3,-5.828){$z_{32}$};
\node (bb3) at (0.7,-5.828){$z_{02}$};
\node (cc3) at (2.7,-5.828){$z_{12}$};
\node (dd3) at (-3.3,-5.828){$z_{22}$};
\fill[black] (a1) circle (4.5pt);
\fill[black] (b1) circle (4.5pt);
\fill[black] (c1) circle (4.5pt);
\fill[black] (d1) circle (4.5pt);

\fill[black] (a2) circle (4.5pt);
\fill[black] (b2) circle (4.5pt);
\fill[black] (c2) circle (4.5pt);
\fill[black] (d2) circle (4.5pt);

\fill[black] (a3) circle (4.5pt);
\fill[black] (b3) circle (4.5pt);
\fill[black] (c3) circle (4.5pt);
\fill[black] (d3) circle (4.5pt);
\path[every node/.style={font=\sffamily\small}]

(a44) edge node{} (a2)
(b44) edge node{} (a3)
(c44) edge node{} (b2)
(d44) edge node{} (b3)
(e44) edge node{} (c2)
(f44) edge node{} (c3)
(g44) edge node{} (d2)
(h44) edge node{} (d3)

(a05) edge [thick,green] node{} (a1)
(b05) edge [thick,green] node{} (b1)
(c05) edge [thick,green] node{} (c1)
(d05) edge [thick,green] node{} (d1)

(a0) edge [thick,red] node{} (a1)
(b0) edge [thick,red] node{} (b1)
(c0) edge [thick,red] node{} (c1)
(d0) edge [thick,red] node{} (d1)

(a1) edge [thick,bend right =30,blue] node{} (b1)
edge [thick,bend right =10,blue] node{} (b1)
edge [thick,bend left =30,blue] node{} (b1)
edge [thick,bend left =10,blue] node{} (b1)
edge [thick, red] node{} (a2)
edge [thick,green] node{} (d2)
edge node{} (a00)

(b1) edge [thick,bend right=30,blue] node{} (c1)
edge [thick,bend right =10,blue] node{} (c1)
edge [thick,bend left =30,blue] node{} (c1)
edge [thick,bend left =10,blue] node{} (c1)
edge [thick,red] node{} (b2)
edge [thick,green] node{} (a2)
edge node{} (b00)
(c1) edge [thick,bend right=0,blue] node{} (c11)
edge [thick,bend right=0,blue] node{} (c12)
edge [thick,bend left=0,blue] node{} (c13)
edge [thick,bend left=0,blue] node{} (c14)
edge [thick,red] node{} (c2)
edge [thick,green] node{} (b2)
edge node{} (c00)
(d11) edge [thick,bend right=0,blue] node{} (d1)
(d12) edge [thick,bend right=0,blue] node{} (d1)
(d13) edge [thick,bend left=0,blue] node{} (d1)
(d14) edge [thick,bend left=0,blue] node{} (d1)
(d1) edge [thick,bend right=30,blue] node{} (a1)
edge [thick,bend right =10,blue] node{} (a1)
edge [thick,bend left =30,blue] node{} (a1)
edge [thick,bend left =10,blue] node{} (a1)
edge [thick,red] node{} (d2)
edge [thick,green] node{} (d6)
edge node{} (d00)

(c6) edge [thick,green] node{} (c2)

(a2) edge [thick,bend right=30,blue] node{} (b2)
edge [thick,bend right =10,blue] node{} (b2)
edge [thick,bend left =30,blue] node{} (b2)
edge [thick,bend left =10,blue] node{} (b2)
edge [thick,red] node{} (a3)
edge [thick,green] node{} (d3)
edge node{} (a01)

(b2) edge [thick,bend right =30,blue] node{} (c2)
edge [thick,bend right =10,blue] node{} (c2)
edge [thick,bend left =30,blue] node{} (c2)
edge [thick,bend left =10,blue] node{} (c2)
edge [thick,red] node{} (b3)
edge [thick,green] node{} (a3)
edge node{} (b01)
(c2)edge [thick,bend right=0,blue] node{} (c21)
edge [thick,bend right=0,blue] node{} (c22)
edge [thick,bend left=0,blue] node{} (c23)
edge [thick,bend left=0,blue] node{} (c24)
edge [thick,red] node{} (c3)
edge [thick,green] node{} (b3)
edge node{} (e00)

(d21) edge [thick,bend right=0,blue] node{} (d2)
(d22) edge [thick,bend right=0,blue] node{} (d2)
(d23) edge [thick,bend left=0,blue] node{} (d2)
(d24) edge [thick,bend left=0,blue] node{} (d2)
(d2) edge [thick,bend right=30,blue] node{} (a2)
edge [thick,bend right =10,blue] node{} (a2)
edge [thick,bend left =30,blue] node{} (a2)
edge [thick,bend left =10,blue] node{} (a2)
edge [thick, red] node{} (d3)
edge [thick,green] node{} (d5)
edge node{} (d01)
(a3) edge [thick,bend right=30,blue] node{} (b3)
edge [thick,bend right =10,blue] node{} (b3)
edge [thick,bend left =30,blue] node{} (b3)
edge [thick,bend left =10,blue] node{} (b3)
edge [thick,red] node{} (a4)
edge node{} (b1)
edge [thick,green] node{} (d55)
(b3) edge [thick,bend right=30,blue] node{} (c3)
edge [thick,bend right =10,blue] node{} (c3)
edge [thick,bend left =30,blue] node{} (c3)
edge [thick,bend left =10,blue] node{} (c3)
edge [thick,red] node{} (b4)
edge node{} (c1)
edge [thick,green] node{} (a55)
(c3) edge [thick,bend right=0,blue] node{} (c31)
edge [thick,bend right=0,blue] node{} (c32)
edge [thick,bend left=0,blue] node{} (c33)
edge [thick,bend left=0,blue] node{} (c34)
edge [thick,red] node{} (c4)
edge node{} (c7)
edge [thick,green] node{} (b55)
(d31) edge [thick,bend right=0,blue] node{} (d3)
(d32) edge [thick,bend right=0,blue] node{} (d3)
(d33) edge [thick,bend left=0,blue] node{} (d3)
(d34) edge [thick,bend left=0,blue] node{} (d3)
(d3) edge [thick,bend right = 30,blue] node{} (a3)
edge [thick,bend right =10,blue] node{} (a3)
edge [thick,bend left =30,blue] node{} (a3)
edge [thick,bend left =10,blue] node{} (a3)
edge [thick,red] node{} (d4)
edge node{} (a1)
edge [thick,green] node{} (d8)
(c5) edge [thick,green] node{} (c3)
(d7) edge node{} (d1);
\end{tikzpicture}$$ Blue lines describe the monodromy action of $v_0,v_1,v_2,v_3$ (which are, in practice, all linear in the horizontal direction), black lines describe the action of $e_0$ while red and green lines describe the action of $e_1,e_2$ respectively.

\bigskip

\hspace{-0.6cm} For a divisor $D \in Div_T(X)$ set $$ \begin{array}{ccc} \vert D \vert_1 :=\sum_{i=0}^s n_i - \sum_{i=0}^r a_i m_i & ; & \vert D \vert_2 = \sum_{i=0}^r m_i \end{array}$$ Set $$ Div^+(k,l) := \left \{ D \vert 0 < k+\vert D \vert_1 \leq s \textrm{ and } 0 < l + \vert D \vert_2 \leq r \right \} \subset Div^+_T(X)$$ For two solutions $(k_1,l_1),(k_2,l_2) \in \mathbb{Z}/(s+1) \mathbb{Z} \oplus \mathbb{Z}/ (r+1) \mathbb{Z}$ we define $$ Hom_{mon}((k_1,l_1),(k_2,l_2) ) := \bigoplus_{D \in M((k_1,l_1),(k_2,l_2)) } \widetilde{M}_D \cdot \mathbb{Z}$$ where $$ M((k_1,l_1),(k_2,l_2)) := \left \{ D \vert \widetilde{M}_D(k_1,l_1) = (k_2,l_2) \textrm{ and } D \in Div^+(k_1,l_1) \right \} $$ Let us note that in terms of these "hom"-functors the M-aligned property could be formulated as follows:

\bigskip

\hspace{-0.6cm} \bf Corollary 4.4 \rm (M-Aligned property (a) II): For any two solutions $(k_1,l_1),(k_2,l_2) \in Crit(X)$ the following holds $$Hom(E_{k_1l_1}, E_{k_2l_2}) \simeq Hom_{mon}((k_1,l_1),(k_2,l_2))$$ Furthermore, the composition map $$ Hom(E_{k_1l_1},E_{k_2l_2}) \otimes Hom(E_{k_2l_2},E_{k_3l_3}) \rightarrow Hom(E_{k_1l_1},E_{k_3l_3})$$ is induced by the map $$Mon((k_1,l_1),(k_2,l_2)) \times Mon((k_2,l_2),(k_3,l_3)) \rightarrow Mon((k_1,l_1),(k_3,l_3)) $$ given by $(D_1,D_2) \mapsto D_1+D_2$.

\bigskip

\hspace{-0.6cm} Motivated by mirror symmetry it is interesting to ask to which extent such "artificially defined" hom-functors between solutions of the asymptotic LG-equations can be related to known geometric structures on the $A$-side, we refer the reader to section 5 for a discussion.

\subsection{M-aligned property for the class (b)} For $r \leq n $ let $ X_{n,r}=Bl_B(\mathbb{P}^{n-r} \times \mathbb{P}^r)$ be the blow up of $\mathbb{P}^{n-r} \times \mathbb{P}^r$ along the multi-linear codimension two subspace $B = \mathbb{P}^{n-r-1} \times \mathbb{P}^{r-1}$. One has $$Div_T(X) = \left ( \bigoplus_{i=1}^n \mathbb{Z} \cdot V_X(e_i) \right )
\bigoplus \left ( \bigoplus_{i=1}^3 \mathbb{Z} \cdot V_X(v_i) \right )$$ Let $L_{nml}=n \cdot V_X(v_1) +m \cdot V_X(v_2) + l \cdot V_X(v_3) \in Pic(X)$ be any element. Then a $T$-divisor satisfies $D \in H^0 (X ; L_{nml})$ if and only if $D \in Div_T^+(X)$ is of the form $$ D = \sum_{i=1}^{n-r}n'_i V_X(e_i) + \sum_{n-r+1}^n m'_i V_X(e_i) + n'' V_X(v_1) + m'' V_X (v_2) + l'' V_X(v_3)$$ where $$ \begin{array}{ccccc} n'' = ( n- \vert n' \vert) & ; & m''=(m-\vert m' \vert) & ; & l'' = (l+ \vert n' \vert + \vert m' \vert +n +m) \end{array}$$

\hspace{-0.6cm} In particular, $\left \{ V(e_1),...,V(e_n),V(v_1),V(v_2),V(v_3) \right \}$ is a basis for the algebra $A_{\mathcal{E}}$. We denote the resulting quiver by $Q(n,r)$. For example, for $X=Bl_B (\mathbb{P}^2 \times \mathbb{P}^2)$ the quiver $Q(4,2)$ is the following

$$
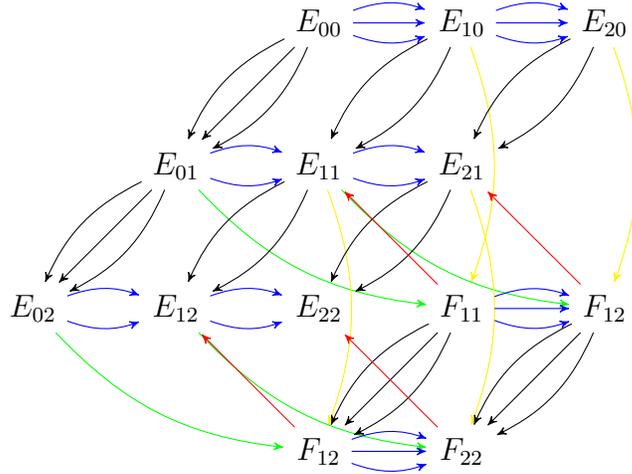
\begin{tikzpicture}[->,>=stealth',shorten >=1pt,auto,node distance=1.9cm ]
\node (1) {$E_{20}$};
\node (4) [left of=1] {$E_{10}$};
\node (7) [left of=4] {$E_{00}$};

\node (2) [below of=4] {$E_{21}$};
\node (5) [left of=2] {$E_{11}$};
\node (8) [left of=5] {$E_{01}$};

\node (3) [below of=5] {$E_{22}$};
\node (6) [left of=3] {$E_{12}$};
\node (9) [left of=6] {$E_{02}$};

\node (10) [right of=3] {$F_{11}$};
\node (11) [right of=10] {$F_{12}$};
\node (12) [below of=10] {$F_{22}$};
\node (13) [left of=12] {$F_{12}$};

\path[every node/.style={font=\sffamily\small}]
(7) edge [bend right=20,blue] node {} (4)
edge [blue] node {} (4)

edge [bend left=20, blue] node {} (4)
edge [bend right=20] node {} (8)
edge node {} (8)
edge [bend left=20] node {} (8)

(8) edge [bend right=20,blue] node {} (5)
edge [bend left=20, blue] node {} (5)

edge [bend right=20] node {} (9)
edge node {} (9)
edge [bend left=20] node {} (9)
edge [bend right=20,green] node{} (10)
(9) edge [bend right=20, blue] node {} (6)

edge [bend left=20, blue] node {} (6)
edge [bend right=20,green] node{} (13)

(4) edge [bend right=20, blue] node {} (1)
edge [ blue] node {} (1)
edge [bend left=20, blue] node {} (1)

edge [bend right=20] node {} (5)

edge [bend left=20] node {} (5)
edge [bend left=20,yellow] node{} (10)

(5) edge [bend right=20, blue] node {} (2)

edge [bend left=20, blue] node {} (2)

edge [bend right=20] node {} (6)

edge [bend left=20] node {} (6)
edge [bend right=20,green] node{} (11)
edge [bend left=20,yellow] node{} (13)

(6) edge [bend right=20, blue] node {} (3)

edge [bend left=20, blue] node {} (3)
edge [bend right=20,green] node{} (12)

(1) edge [bend right=20] node {} (2)

edge [bend left=20] node {} (2)
edge [bend left=20,yellow] node{} (11)

(2) edge [bend right=20] node {} (3)

edge [bend left=20] node {} (3)
edge [bend left=20,yellow] node{} (12)

(10)
edge [bend right=20, blue] node {} (11)
edge [blue] node {} (11)
edge [bend left=20, blue] node {} (11)
edge [red] node {} (5)
edge [ bend right=20] node {} (13)
edge node {} (13)
edge [bend left=20] node {} (13)

(11) edge [red] node {} (2)
edge [bend right=20] node {} (12)
edge node {} (12)
edge [bend left=20] node {} (12)

(12) edge [red] node {} (3)

(13)
edge [bend right=20, blue] node {} (12)
edge [ blue] node {} (12)
edge [bend left=20, blue] node {} (12)
edge [red] node {} (6);

\end{tikzpicture}
$$

\hspace{-0.6cm} Where two of the blue lines stand for the elements $V_X(e_1), V_X(e_2)$, two of the black lines for $V_X(e_3),V_X(e_4)$, green lines for $V_X(v_1)$, yellow for $V_X(v_2)$ and red for $V_X(v_3)$. The additional blue and black lines stand for $V_X(v_1) + V_X(v_3)$ and $V_X(v_2) +V_X(v_3)$ respectively.

\hspace{-0.6cm} On the other hand, on the Landau-Ginzburg side, the base point is given by the system of equations $$ W^t_X : = \left \{ \begin{array}{c} e^{-t} z_i - \frac{1}{ \prod_{i=1}^{n-r} z_i } - \frac{1}{\prod_{i=1}^{n-r} z_i \cdot \prod_{i=1}^r w_i } =0 \\ e^{-t} w_j - \frac{1}{ \prod_{i=1}^{r} w_i } - \frac{1}{\prod_{i=1}^{n-r} z_i \cdot \prod_{i=1}^r w_i } =0 \end{array} \right \}_{i=1,j=1}^{n-r,r} $$ For $ 1 \leq i \leq n-r$, $n-r+1 \leq j \leq n$ and $ \theta \in [0,1)$ we define the following loops in $Gr(X) \setminus R_X$:

$$ W^t_{V_X(e_i)}( \theta) : = \left \{ \begin{array}{c} e^{2 \pi i \delta(i,i') \theta} z_{i'} - \frac{1}{ \prod_{i=1}^{n-r} z_i } - \frac{1}{\prod_{i=1}^{n-r} z_i \cdot \prod_{i=1}^r w_i } =0 \\ e^{-t} w_j - \frac{1}{ \prod_{i=1}^{r} w_i } - \frac{1}{\prod_{i'=1}^{n-r} z_i \cdot \prod_{i=1}^r w_i } =0 \end{array} \right \}_{i'=1,j=1}^{n-r,r} $$

$$ W^t_{V_X(e_j)}( \theta) : = \left \{ \begin{array}{c} e^{-t} z_i - \frac{1}{ \prod_{i=1}^{n-r} z_i } - \frac{1}{\prod_{i=1}^{n-r} z_i \cdot \prod_{i=1}^r w_i } =0 \\ e^{2 \pi i \delta(j,j') \theta} w_{j'} - \frac{1}{ \prod_{i=1}^{r} w_i } - \frac{1}{\prod_{i=1}^{n-r} z_i \cdot \prod_{i=1}^r w_i } =0 \end{array} \right \}_{i=1,j'=1}^{n-r,r} $$

$$ \widetilde{W}^t_{V_X(e_i)-V_X(v_3)}( \theta) : = \left \{ \begin{array}{c} e^{2 \pi i \delta(i,i') \theta} z_{i'} - \frac{1}{ \prod_{i=1}^{n-r} z_i } - \frac{1}{\prod_{i=1}^{n-r} z_i \cdot \prod_{i=1}^r w_i } =0 \\ e^{-t} w_j - \frac{1}{ \prod_{i=1}^{r} w_i } - \frac{e^{-2 \pi i \theta}}{\prod_{i'=1}^{n-r} z_i \cdot \prod_{i=1}^r w_i } =0 \end{array} \right \}_{i'=1,j=1}^{n-r,r} $$

$$ \widetilde{W}^t_{V_X(e_j)-V_X(v_3)}( \theta) : = \left \{ \begin{array}{c} e^{-t} z_i - \frac{1}{ \prod_{i=1}^{n-r} z_i } - \frac{e^{-2 \pi i \theta}}{\prod_{i=1}^{n-r} z_i \cdot \prod_{i=1}^r w_i } =0 \\ e^{2 \pi i \delta(j,j') \theta} w_{j'} - \frac{1}{ \prod_{i=1}^{r} w_i } - \frac{1}{\prod_{i'=1}^{n-r} z_i \cdot \prod_{i=1}^r w_i } =0 \end{array} \right \}_{i=1,j'=1}^{n-r,r} $$

$$ \widetilde{W}^t_{V_X(v_3)}( \theta) : = \left \{ \begin{array}{c} z_i - \frac{1}{ \prod_{i=1}^{n-r} z_i } - \frac{1}{\prod_{i=1}^{n-r} z_i \cdot \prod_{i=1}^r w_i } =0 \\ e^{-t} w_j - \frac{1}{ \prod_{i=1}^{r} w_i } - \frac{e^{2 \pi i \theta}}{\prod_{i'=1}^{n-r} z_i \cdot \prod_{i=1}^r w_i } =0 \end{array} \right \}_{i=1,j=1}^{n-r,r} $$

$$ \widetilde{W'}^t_{V_X(v_3)}( \theta) : = \left \{ \begin{array}{c} e^{-t} z_{i'} - \frac{1}{ \prod_{i=1}^{n-r} z_i } - \frac{e^{-2 \pi i \theta}}{\prod_{i=1}^{n-r} z_i \cdot \prod_{i=1}^r w_i } =0 \\ w_j - \frac{1}{ \prod_{i=1}^{r} w_i } - \frac{1}{\prod_{i'=1}^{n-r} z_i \cdot \prod_{i=1}^r w_i } =0 \end{array} \right \}_{i=1,j=1}^{n-r,r} $$ Set also $$ U^t(\theta) : = \left \{ \begin{array}{c} e^{-t} z_i - \frac{1}{ \prod_{i=1}^{n-r} z_i } - \frac{1}{\prod_{i=1}^{n-r} z_i \cdot \prod_{i=1}^r w_i } =0 \\ e^{-t (1- \theta)} w_j - \frac{1}{ \prod_{i=1}^{r} w_i } - \frac{1}{\prod_{i=1}^{n-r} z_i \cdot \prod_{i=1}^r w_i } =0 \end{array} \right \}_{i=1,j=1}^{n-r,r} $$ $$ V^t( \theta) : = \left \{ \begin{array}{c} e^{-t (1- \theta)} z_i - \frac{1}{ \prod_{i=1}^{n-r} z_i } - \frac{1}{\prod_{i=1}^{n-r} z_i \cdot \prod_{i=1}^r w_i } =0 \\ e^{-t} w_j - \frac{1}{ \prod_{i=1}^{r} w_i } - \frac{1}{\prod_{i=1}^{n-r} z_i \cdot \prod_{i=1}^r w_i } =0 \end{array} \right \}_{i=1,j=1}^{n-r,r}$$

\hspace{-0.6cm} Set $$\begin{array}{cccc} \Gamma_{V_X(e_i)}: = [(U^t)^{-1} \circ W^t_{V_X(e_i)} \circ U^t] & ; & \Gamma_{V_X(e_j)}:=[(V^t)^{-1} \circ W^t_{V_X(e_j)} \circ V^t] & ; \end{array} $$
$$\begin{array}{cccc} \Gamma_{V_X(v_1)}:= [(U^t)^{-1} \circ \widetilde{W}^t_{V_X(e_i)-V_X(v_3)} \circ U^t] & ; & \Gamma_{V_X(v_2)} :=
[(V^t)^{-1} \circ \widetilde{W}^t_{V_X(e_j)-V_X(v_3)} \circ V^t] & ; \end{array} $$
$$\begin{array}{ccc} \Gamma_{V_X(v_3)}:= [(U^t)^{-1} \circ \widetilde{W}^t_{V_X(v_3)} \circ U^t] & ; & \Gamma'_{V_X(v_3)} :=
[(V^t)^{-1} \circ \widetilde{W'}^t_{V_X(v_3)} \circ V^t] \end{array} $$ Let $$ Crit(X_{n,r} ; W^t_X) = \left \{ (z_k,w_l) \right \}_{k=0,l=0}^{n-r,r} \bigcup \left \{ (z'_m,w'_n) \right \}_{m=0,n=0}^{n-r-1,r-1} $$ be as in Definition 3.2.1. The following example illustrates the $M$-aligned property for the class (b):

\bigskip

\hspace{-0.6cm} \bf Example 4.5 \rm (monodromies for $X= Bl_B(\mathbb{P}^2 \times \mathbb{P}^2)$ with $B= \mathbb{P}^1 \times \mathbb{P}^1$): The exceptional map of Definition 3.2.1 is defined in terms of the solutions of the LG-system for $$ f_t(z,w) = e^{-t} z_1+e^{-t}z_2+e^{-t}w_1+e^{-t}w_2 + \frac{1}{z_1z_2} + \frac{1}{w_1 w_2 } + \frac{1}{z_1z_2w_1w_2} $$ for $0<<t$. The following is a graphical description of the $z_1$ and $w_1$ coordinates of the solution scheme $Crit(X ; W_X^t)$:

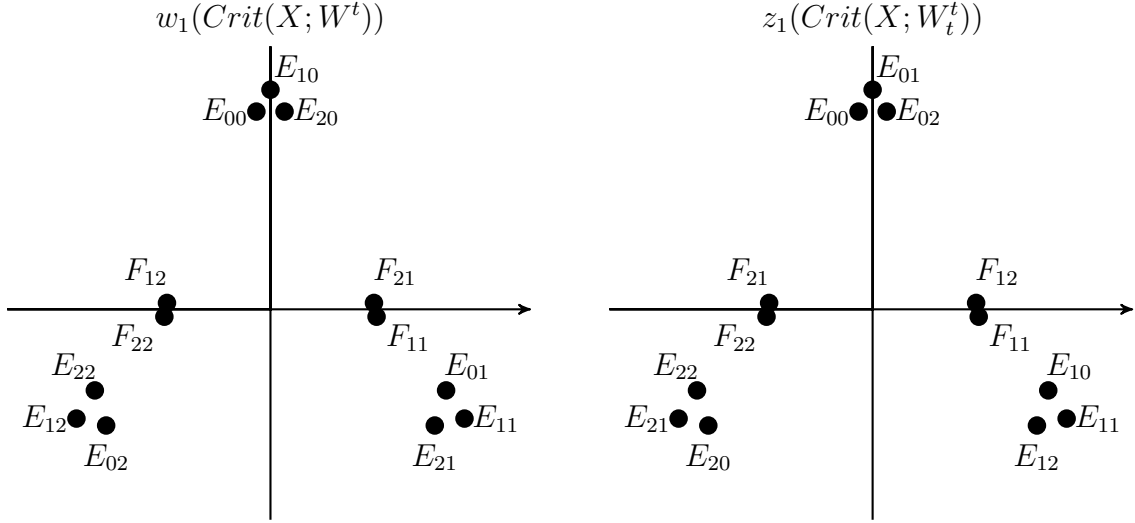
$$ \begin{array}{ccc} \begin{tikzpicture}[->,>=stealth',shorten >=1pt,auto,node distance=3.5cm,main node/.style={circle,draw,font=\sffamily\bfseries \small}, scale=1.4]
\draw [->,thick] (0,-2)--(0,2.5) node (yaxis) [above] {$w_1(Crit(X; W^t))$}
|- (-2.5,0)--(2.5,0) node (xaxis) [right] {};

\node (a) at (0, 2.08669){};
\node (b) at (1.56015, -1.10404){};
\node (c) at (-1.56015, -1.10404){};
\node (d) at (1.84216, -1.03819){};
\node (e) at (-1.84216, -1.03819){};
\node (f) at (0.133691, 1.87899){};
\node (g) at (-0.133691, 1.87899){};
\node (h) at (1.66754, -0.770428){};
\node (h1) at (-1.66754, -0.770428){};
\node (h2) at (0.982687, 0.0606989){};
\node (h3) at (-0.982687, 0.0606989){};
\node (h4) at (1.0074, -0.0703753){};
\node (h5) at (-1.0074, -0.0703753){};

\fill[black] (a) circle (2.5pt);
\fill[black] (b) circle (2.5pt);
\fill[black] (c) circle (2.5pt);
\fill[black] (d) circle (2.5pt);
\fill[black] (e) circle (2.5pt);
\fill[black] (f) circle (2.5pt);
\fill[black] (g) circle (2.5pt);
\fill[black] (h) circle (2.5pt);
\fill[black] (h1) circle (2.5pt);
\fill[black] (h2) circle (2.5pt);
\fill[black] (h3) circle (2.5pt);
\fill[black] (h4) circle (2.5pt);
\fill[black] (h5) circle (2.5pt);

\node (a1) at (0.25, 2.28669){$E_{10}$};
\node (b1) at (1.56015, -1.40404){$E_{21}$};
\node (c1) at (-1.56015, -1.40404){$E_{02}$};
\node (d1) at (2.14216, -1.03819){$E_{11}$};
\node (e1) at (-2.14216, -1.03819){$E_{12}$};
\node (f1) at (0.433691, 1.87899){$E_{20}$};
\node (g1) at (-0.433691, 1.87899){$E_{00}$};
\node (h1) at (1.86754, -0.540428){$E_{01}$};
\node (h11) at (-1.86754, -0.540428){$E_{22}$};
\node (h21) at (1.182687, 0.3606989){$F_{21}$};
\node (h31) at (-1.182687, 0.3606989){$F_{12}$};
\node (h41) at (1.3074, -0.2703753){$F_{11}$};
\node (h51) at (-1.3074, -0.2703753){$F_{22}$};

\end{tikzpicture}
& &
\begin{tikzpicture}[->,>=stealth',shorten >=1pt,auto,node distance=3.5cm,main node/.style={circle,draw,font=\sffamily\bfseries \small}, scale=1.4]
\draw [->,thick] (0,-2)--(0,2.5) node (yaxis) [above] {$z_1(Crit(X;W^t_t))$}
|- (-2.5,0)--(2.5,0) node (xaxis) [right] {};

\node (a) at (0, 2.08669){};
\node (b) at (1.56015, -1.10404){};
\node (c) at (-1.56015, -1.10404){};
\node (d) at (1.84216, -1.03819){};
\node (e) at (-1.84216, -1.03819){};
\node (f) at (0.133691, 1.87899){};
\node (g) at (-0.133691, 1.87899){};
\node (h) at (1.66754, -0.770428){};
\node (h1) at (-1.66754, -0.770428){};
\node (h2) at (0.982687, 0.0606989){};
\node (h3) at (-0.982687, 0.0606989){};
\node (h4) at (1.0074, -0.0703753){};
\node (h5) at (-1.0074, -0.0703753){};

\fill[black] (a) circle (2.5pt);
\fill[black] (b) circle (2.5pt);
\fill[black] (c) circle (2.5pt);
\fill[black] (d) circle (2.5pt);
\fill[black] (e) circle (2.5pt);
\fill[black] (f) circle (2.5pt);
\fill[black] (g) circle (2.5pt);
\fill[black] (h) circle (2.5pt);
\fill[black] (h1) circle (2.5pt);
\fill[black] (h2) circle (2.5pt);
\fill[black] (h3) circle (2.5pt);
\fill[black] (h4) circle (2.5pt);
\fill[black] (h5) circle (2.5pt);

\node (a1) at (0.25, 2.28669){$E_{01}$};
\node (b1) at (1.56015, -1.40404){$E_{12}$};
\node (c1) at (-1.56015, -1.40404){$E_{20}$};
\node (d1) at (2.14216, -1.03819){$E_{11}$};
\node (e1) at (-2.14216, -1.03819){$E_{21}$};
\node (f1) at (0.433691, 1.87899){$E_{02}$};
\node (g1) at (-0.433691, 1.87899){$E_{00}$};
\node (h1) at (1.86754, -0.540428){$E_{10}$};
\node (h11) at (-1.86754, -0.540428){$E_{22}$};
\node (h21) at (1.182687, 0.3606989){$F_{12}$};
\node (h31) at (-1.182687, 0.3606989){$F_{21}$};
\node (h41) at (1.3074, -0.2703753){$F_{11}$};
\node (h51) at (-1.3074, -0.2703753){$F_{22}$};

\end{tikzpicture}
\end{array}
$$ The action of $U^t(\theta)$ and $V^t(\theta)$ for $[0,1)$ is described as follows:

$$ \begin{array}{ccc} 
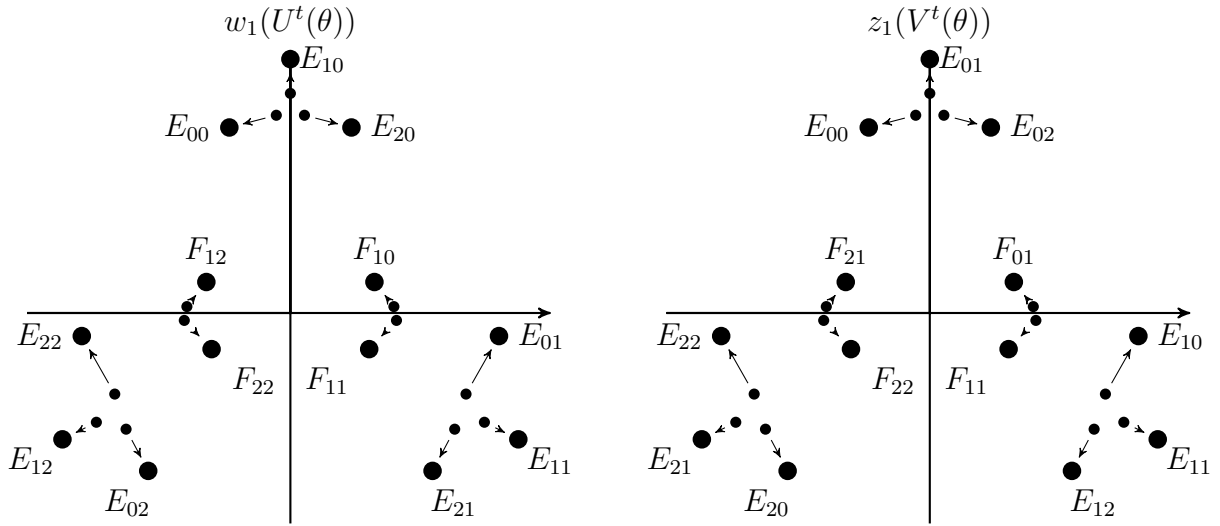
\begin{tikzpicture}[->,>=stealth',shorten >=1pt,auto,node distance=3.5cm,main node/.style={circle,draw,font=\sffamily\bfseries \small}, scale=1.4]
\draw [->,thick] (0,-2)--(0,2.5) node (yaxis) [above] {$w_1(U^t( \theta))$}
|- (-2.5,0)--(2.5,0) node (xaxis) [right] {};

\node (a) at (0, 2.08669){};
\node (b) at (1.56015, -1.10404){};
\node (c) at (-1.56015, -1.10404){};
\node (d) at (1.84216, -1.03819){};
\node (e) at (-1.84216, -1.03819){};
\node (f) at (0.133691, 1.87899){};
\node (g) at (-0.133691, 1.87899){};
\node (h) at (1.66754, -0.770428){};
\node (h1) at (-1.66754, -0.770428){};
\node (h2) at (0.982687, 0.0606989){};
\node (h3) at (-0.982687, 0.0606989){};
\node (h4) at (1.0074, -0.0703753){};
\node (h5) at (-1.0074, -0.0703753){};

\node (a1) at (0, 2.41225){};
\node (b1) at (1.34988, -1.50026){};
\node (c1) at (-1.34988, -1.50026){};
\node (d1) at (0.797981, 0.294386){};
\node (e1) at (-0.797981, 0.294386){};
\node (f1) at (0.579635, 1.76278){};
\node (g1) at (-0.579635, 1.76278){};
\node (h01) at (2.16423, -1.20049){};
\node (h11) at (-2.16423, -1.20049){};
\node (h21) at (-0.748239, -0.343339){};
\node (h31) at (0.748239, -0.343339){};
\node (h41) at (-1.98059, -0.219197){};
\node (h51) at (1.98059, -0.219197){};

\node (a11) at (0.3, 2.41225){$E_{10}$};
\node (b11) at (1.54988, -1.80026){$E_{21}$};
\node (c11) at (-1.54988, -1.80026){$E_{02}$};
\node (d11) at (0.797981, 0.594386){$F_{10}$};
\node (e11) at (-0.797981, 0.594386){$F_{12}$};
\node (f11) at (0.979635, 1.76278){$E_{20}$};
\node (g11) at (-0.979635, 1.76278){$E_{00}$};
\node (h011) at (2.46423, -1.40049){$E_{11}$};
\node (h111) at (-2.46423, -1.40049){$E_{12}$};
\node (h211) at (-0.348239, -0.643339){$F_{22}$};
\node (h311) at (0.348239, -0.643339){$F_{11}$};
\node (h411) at (-2.38059, -0.219197){$E_{22}$};
\node (h511) at (2.38059, -0.219197){$E_{01}$};

\fill[black] (a) circle (1.5pt);
\fill[black] (b) circle (1.5pt);
\fill[black] (c) circle (1.5pt);
\fill[black] (d) circle (1.5pt);
\fill[black] (e) circle (1.5pt);
\fill[black] (f) circle (1.5pt);
\fill[black] (g) circle (1.5pt);
\fill[black] (h) circle (1.5pt);
\fill[black] (h1) circle (1.5pt);
\fill[black] (h2) circle (1.5pt);
\fill[black] (h3) circle (1.5pt);
\fill[black] (h4) circle (1.5pt);
\fill[black] (h5) circle (1.5pt);

\fill[black] (a1) circle (2.5pt);
\fill[black] (b1) circle (2.5pt);
\fill[black] (c1) circle (2.5pt);
\fill[black] (d1) circle (2.5pt);
\fill[black] (e1) circle (2.5pt);
\fill[black] (f1) circle (2.5pt);
\fill[black] (g1) circle (2.5pt);
\fill[black] (h01) circle (2.5pt);
\fill[black] (h11) circle (2.5pt);
\fill[black] (h21) circle (2.5pt);
\fill[black] (h31) circle (2.5pt);
\fill[black] (h41) circle (2.5pt);
\fill[black] (h51) circle (2.5pt);

\path[every node/.style={font=\sffamily\small}]
(a) edge (a1)
(b) edge (b1)
(c) edge (c1)
(d) edge (h01)
(e) edge (h11)
(f) edge (f1)
(g) edge (g1)
(h) edge (h51)
(h1) edge (h41)
(h2) edge (d1)
(h3) edge (e1)
(h4) edge (h31)
(h5) edge (h21);
\end{tikzpicture}
& &
\begin{tikzpicture}[->,>=stealth',shorten >=1pt,auto,node distance=3.5cm,main node/.style={circle,draw,font=\sffamily\bfseries \small}, scale=1.4]
\draw [->,thick] (0,-2)--(0,2.5) node (yaxis) [above] {$z_1(V^t(\theta))$}
|- (-2.5,0)--(2.5,0) node (xaxis) [right] {};

\node (a) at (0, 2.08669){};
\node (b) at (1.56015, -1.10404){};
\node (c) at (-1.56015, -1.10404){};
\node (d) at (1.84216, -1.03819){};
\node (e) at (-1.84216, -1.03819){};
\node (f) at (0.133691, 1.87899){};
\node (g) at (-0.133691, 1.87899){};
\node (h) at (1.66754, -0.770428){};
\node (h1) at (-1.66754, -0.770428){};
\node (h2) at (0.982687, 0.0606989){};
\node (h3) at (-0.982687, 0.0606989){};
\node (h4) at (1.0074, -0.0703753){};
\node (h5) at (-1.0074, -0.0703753){};

\node (a1) at (0, 2.41225){};
\node (b1) at (1.34988, -1.50026){};
\node (c1) at (-1.34988, -1.50026){};
\node (d1) at (0.797981, 0.294386){};
\node (e1) at (-0.797981, 0.294386){};
\node (f1) at (0.579635, 1.76278){};
\node (g1) at (-0.579635, 1.76278){};
\node (h01) at (2.16423, -1.20049){};
\node (h11) at (-2.16423, -1.20049){};
\node (h21) at (-0.748239, -0.343339){};
\node (h31) at (0.748239, -0.343339){};
\node (h41) at (-1.98059, -0.219197){};
\node (h51) at (1.98059, -0.219197){};

\node (a11) at (0.3, 2.41225){$E_{01}$};
\node (b11) at (1.54988, -1.80026){$E_{12}$};
\node (c11) at (-1.54988, -1.80026){$E_{20}$};
\node (d11) at (0.797981, 0.594386){$F_{01}$};
\node (e11) at (-0.797981, 0.594386){$F_{21}$};
\node (f11) at (0.979635, 1.76278){$E_{02}$};
\node (g11) at (-0.979635, 1.76278){$E_{00}$};
\node (h011) at (2.46423, -1.40049){$E_{11}$};
\node (h111) at (-2.46423, -1.40049){$E_{21}$};
\node (h211) at (-0.348239, -0.643339){$F_{22}$};
\node (h311) at (0.348239, -0.643339){$F_{11}$};
\node (h411) at (-2.38059, -0.219197){$E_{22}$};
\node (h511) at (2.38059, -0.219197){$E_{10}$};

\fill[black] (a) circle (1.5pt);
\fill[black] (b) circle (1.5pt);
\fill[black] (c) circle (1.5pt);
\fill[black] (d) circle (1.5pt);
\fill[black] (e) circle (1.5pt);
\fill[black] (f) circle (1.5pt);
\fill[black] (g) circle (1.5pt);
\fill[black] (h) circle (1.5pt);
\fill[black] (h1) circle (1.5pt);
\fill[black] (h2) circle (1.5pt);
\fill[black] (h3) circle (1.5pt);
\fill[black] (h4) circle (1.5pt);
\fill[black] (h5) circle (1.5pt);

\fill[black] (a1) circle (2.5pt);
\fill[black] (b1) circle (2.5pt);
\fill[black] (c1) circle (2.5pt);
\fill[black] (d1) circle (2.5pt);
\fill[black] (e1) circle (2.5pt);
\fill[black] (f1) circle (2.5pt);
\fill[black] (g1) circle (2.5pt);
\fill[black] (h01) circle (2.5pt);
\fill[black] (h11) circle (2.5pt);
\fill[black] (h21) circle (2.5pt);
\fill[black] (h31) circle (2.5pt);
\fill[black] (h41) circle (2.5pt);
\fill[black] (h51) circle (2.5pt);

\path[every node/.style={font=\sffamily\small}]
(a) edge (a1)
(b) edge (b1)
(c) edge (c1)
(d) edge (h01)
(e) edge (h11)
(f) edge (f1)
(g) edge (g1)
(h) edge (h51)
(h1) edge (h41)
(h2) edge (d1)
(h3) edge (e1)
(h4) edge (h31)
(h5) edge (h21);

\end{tikzpicture}
\end{array}
$$

\hspace{-0.6cm} The following is a graphical description of the corresponding monodromy actions:

$$ \begin{array}{ccc} 
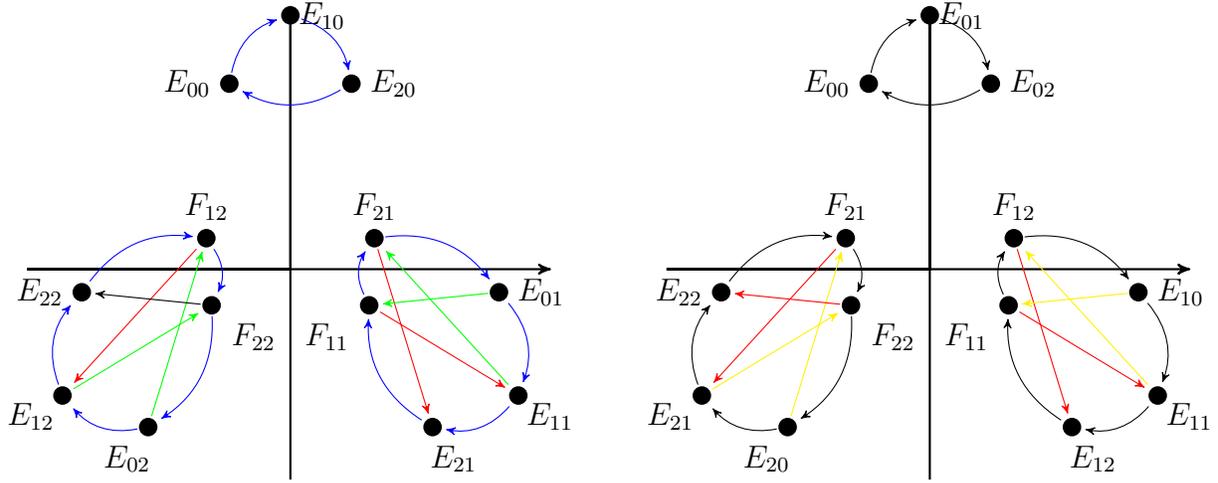
\begin{tikzpicture}[->,>=stealth',shorten >=1pt,auto,node distance=3.5cm,main node/.style={circle,draw,font=\sffamily\bfseries \small}, scale=1.4]
\draw [->,thick] (0,-2)--(0,2.5) node {}(yaxis) [above]
|- (-2.5,0)--(2.5,0) node {}(xaxis) [right] {} ;

\node (a1) at (0, 2.41225){};
\node (b1) at (1.34988, -1.50026){};
\node (c1) at (-1.34988, -1.50026){};
\node (d1) at (0.797981, 0.294386){};
\node (e1) at (-0.797981, 0.294386){};
\node (f1) at (0.579635, 1.76278){};
\node (g1) at (-0.579635, 1.76278){};
\node (h01) at (2.16423, -1.20049){};
\node (h11) at (-2.16423, -1.20049){};
\node (h21) at (-0.748239, -0.343339){};
\node (h31) at (0.748239, -0.343339){};
\node (h41) at (-1.98059, -0.219197){};
\node (h51) at (1.98059, -0.219197){};

\node (a11) at (0.3, 2.41225){$E_{10}$};
\node (b11) at (1.54988, -1.80026){$E_{21}$};
\node (c11) at (-1.54988, -1.80026){$E_{02}$};
\node (d11) at (0.797981, 0.594386){$F_{21}$};
\node (e11) at (-0.797981, 0.594386){$F_{12}$};
\node (f11) at (0.979635, 1.76278){$E_{20}$};
\node (g11) at (-0.979635, 1.76278){$E_{00}$};
\node (h011) at (2.46423, -1.40049){$E_{11}$};
\node (h111) at (-2.46423, -1.40049){$E_{12}$};
\node (h211) at (-0.348239, -0.643339){$F_{22}$};
\node (h311) at (0.348239, -0.643339){$F_{11}$};
\node (h411) at (-2.38059, -0.219197){$E_{22}$};
\node (h511) at (2.38059, -0.219197){$E_{01}$};

\fill[black] (a1) circle (2.5pt);
\fill[black] (b1) circle (2.5pt);
\fill[black] (c1) circle (2.5pt);
\fill[black] (d1) circle (2.5pt);
\fill[black] (e1) circle (2.5pt);
\fill[black] (f1) circle (2.5pt);
\fill[black] (g1) circle (2.5pt);
\fill[black] (h01) circle (2.5pt);
\fill[black] (h11) circle (2.5pt);
\fill[black] (h21) circle (2.5pt);
\fill[black] (h31) circle (2.5pt);
\fill[black] (h41) circle (2.5pt);
\fill[black] (h51) circle (2.5pt);

\path[every node/.style={font=\sffamily\small}]
(a1) edge [blue,->,bend left] (f1)
(b1) edge [blue,->,bend left] (h31)
(c1) edge [blue,bend left] (h11)
edge [green] (e1)
(d1) edge [blue,->,bend left] (h51)
edge [red] (b1)
(e1) edge [blue,->,bend left] (h21)
edge [red] (h11)
(f1) edge [blue,->,bend left] (g1)
(g1) edge [blue,->,bend left] (a1)
(h01) edge [blue,->,bend left] (b1)
edge [green] (d1)
(h11) edge [blue,bend left] (h41)
edge [green] (h21)
(h21) edge [blue,->,bend left] (c1)
edge (h41)
(h31) edge [blue,->,bend left] (d1)
edge [red] (h01)
(h41) edge [blue,->,bend left] (e1)
(h51) edge [blue,->,bend left] (h01)
edge [green] (h31);
\end{tikzpicture}
& &
\begin{tikzpicture}[->,>=stealth',shorten >=1pt,auto,node distance=3.5cm,main node/.style={circle,draw,font=\sffamily\bfseries \small}, scale=1.4]
\draw [->,thick] (0,-2)--(0,2.5) node {} (yaxis) [above]
|- (-2.5,0)--(2.5,0) node {} (xaxis) [right] {};

\node (a1) at (0, 2.41225){};
\node (b1) at (1.34988, -1.50026){};
\node (c1) at (-1.34988, -1.50026){};
\node (d1) at (0.797981, 0.294386){};
\node (e1) at (-0.797981, 0.294386){};
\node (f1) at (0.579635, 1.76278){};
\node (g1) at (-0.579635, 1.76278){};
\node (h01) at (2.16423, -1.20049){};
\node (h11) at (-2.16423, -1.20049){};
\node (h21) at (-0.748239, -0.343339){};
\node (h31) at (0.748239, -0.343339){};
\node (h41) at (-1.98059, -0.219197){};
\node (h51) at (1.98059, -0.219197){};

\node (a11) at (0.3, 2.41225){$E_{01}$};
\node (b11) at (1.54988, -1.80026){$E_{12}$};
\node (c11) at (-1.54988, -1.80026){$E_{20}$};
\node (d11) at (0.797981, 0.594386){$F_{12}$};
\node (e11) at (-0.797981, 0.594386){$F_{21}$};
\node (f11) at (0.979635, 1.76278){$E_{02}$};
\node (g11) at (-0.979635, 1.76278){$E_{00}$};
\node (h011) at (2.46423, -1.40049){$E_{11}$};
\node (h111) at (-2.46423, -1.40049){$E_{21}$};
\node (h211) at (-0.348239, -0.643339){$F_{22}$};
\node (h311) at (0.348239, -0.643339){$F_{11}$};
\node (h411) at (-2.38059, -0.219197){$E_{22}$};
\node (h511) at (2.38059, -0.219197){$E_{10}$};

\fill[black] (a1) circle (2.5pt);
\fill[black] (b1) circle (2.5pt);
\fill[black] (c1) circle (2.5pt);
\fill[black] (d1) circle (2.5pt);
\fill[black] (e1) circle (2.5pt);
\fill[black] (f1) circle (2.5pt);
\fill[black] (g1) circle (2.5pt);
\fill[black] (h01) circle (2.5pt);
\fill[black] (h11) circle (2.5pt);
\fill[black] (h21) circle (2.5pt);
\fill[black] (h31) circle (2.5pt);
\fill[black] (h41) circle (2.5pt);
\fill[black] (h51) circle (2.5pt);

\path[every node/.style={font=\sffamily\small}]
(a1) edge [->,bend left] (f1)
(b1) edge [->,bend left] (h31)
(c1) edge [bend left] (h11)
edge [yellow] (e1)
(d1) edge [->,bend left] (h51)
edge [red] (b1)
(e1) edge [->,bend left] (h21)
edge [red] (h11)
(f1) edge [->,bend left] (g1)
(g1) edge [->,bend left] (a1)
(h01) edge [->,bend left] (b1)
edge [yellow] (d1)
(h11) edge [bend left] (h41)
edge [yellow] (h21)
(h21) edge [->,bend left] (c1)
edge [red] (h41)
(h31) edge [->,bend left] (d1)
edge [red] (h01)
(h41) edge [->,bend left] (e1)
(h51) edge [->,bend left] (h01)
edge [yellow] (h31);
\end{tikzpicture}
\end{array}
$$

\hspace{-0.6cm} Blue lines represent the monodromies of $[\widetilde{W}^t_{V_X(e_i)}]$ for $i=1,...,n-r$, black lines $[\widetilde{W}^t_{V_X(e_j)}]$ for $j=n-r+1,...,n$, green lines represent the monodromies
of $[\widetilde{W}^t_{V_X(v_1)}]$, yellow lines $[\widetilde{W}^t_{V_X(v_2)}]$ and red lines the monodromies of $[\widetilde{W}^t_{V_X(v_3)}]=[\widetilde{W}'^t_{V_X(v_3)}]$. In particular, one has

\bigskip

(1) $M( \Gamma_{V_X(e_i)}) (z_k,w_l) = (z_{k+1},w_l) $ and $M(\Gamma_{V_X(e_i)} ) (z'_m,w'_n)= (z'_{m+1},w'_n)$.

\bigskip

(2) $M( \Gamma_{V_X(e_j)}) (z_k,w_l) = (z_k,w_{l+1}) $ and $M(\Gamma_{V_X(e_j)} ) (z'_m,w'_n)= (z'_m,w'_{n+1})$.

\bigskip

(3) $M(\Gamma_{V_X(v_1)}) (z_k,w_l) = (z'_{k+1},w'_l)$.

\bigskip

(4) $M(\Gamma_{V_X(v_2)}) (z_k,w_l) = (z'_k,w'_{l+1})$.

\bigskip

(5) $M(\Gamma_{V_X(v_3)}) (z'_k,w'_l) = M(\Gamma'_{V_X(v_3)})(z'_k,w'_l)=(z_k,w_l)$.

\bigskip

\hspace{-0.6cm} Relations (1)-(5) establish the $M$-aligned property in this case.

\subsection{M-aligned property for the class (c)} For $b < n-1 $ let $ X_{n,b}=Bl_B(
\mathbb{P}(\mathcal{O}_{\mathbb{P}^{n-1}}) \oplus \mathcal{O}_{\mathbb{P}^{n-1}}(b))$ be the blow up of the projective bundle $\mathbb{P}(\mathcal{O}_{\mathbb{P}^{n-1}} \oplus \mathcal{O}_{\mathbb{P}^{n-1}}(b))$ along the codimension two subspace $B = \mathbb{P}^{n-2}$. One has $$Div_T(X) = \left ( \bigoplus_{i=1}^n \mathbb{Z} \cdot V_X(e_i) \right )
\bigoplus \left ( \bigoplus_{i=1}^3 \mathbb{Z} \cdot V_X(u_i) \right )$$ Let $L_{nml}=n \cdot V +m T + l \cdot Y \in Pic(X)$ be any element. Then a $T$-divisor satisfies $D \in H^0 (X ; L_{nml})$ if and only if $D \in Div_T^+(X)$ is of the form $$ D = \sum_{i=1}^{n-1} n'_i V_X(e_i) + m' V_X(e_n) + n'' V_X(u_1) + m'' V_X (u_2) + l'' V_X(u_3)$$ where $$ \begin{array}{ccccc} n = ( \vert n' \vert + m' +m'' \cdot b) & ; & m=(m''-m'+l'') & ; & l = (n''+m'') \end{array}$$

\hspace{-0.6cm} In particular, $\left \{ V(e_1),...,V(e_n),V(u_1),V(u_2),V(u_3) \right \}$ is a basis for the algebra $A_{\mathcal{E}}$. We denote the resulting quiver by $Q(n,r)$. For example, for $X=Bl_B (\mathbb{P}^2 \times \mathbb{P}^2)$ the quiver $Q(4,2)$ is the following
$$
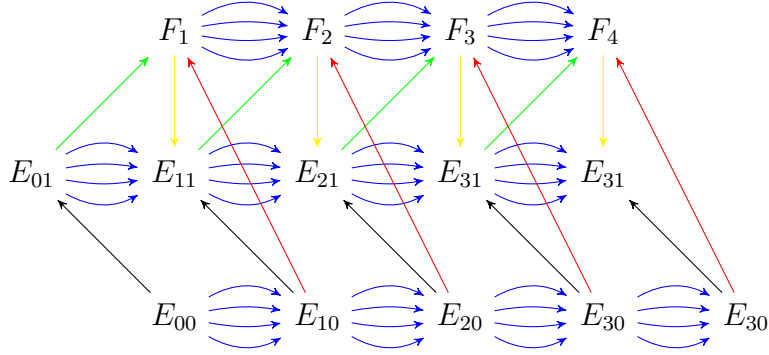
\begin{tikzpicture}[->,>=stealth',shorten >=1pt,auto,node distance=1.9cm]
\node (1) {$E_{01}$};
\node (4) [right of=1] {$E_{11}$};
\node (7) [right of=4] {$E_{21}$};
\node (10) [right of=7] {$E_{31}$};
\node (13) [right of=10] {$E_{31}$};

\node (2) [below of=4] {$E_{00}$};
\node (5) [right of=2] {$E_{10}$};
\node (8) [right of=5] {$E_{20}$};
\node (11) [right of=8] {$E_{30}$};
\node (14) [right of=11] {$E_{30}$};

\node (3) [above of =4] {$F_1$};
\node (6) [above of =7] {$F_2$};
\node (9) [above of=10] {$F_3$};
\node (12) [above of=13] {$F_4$};

\path[every node/.style={font=\sffamily\small}]
(1) edge [bend right=30,blue] node {} (4)
edge [bend right =10,blue] node {} (4)
edge [bend left=30, blue] node {} (4)
edge [bend left = 10,blue] node {} (4)
edge [green] node {} (3)

(2) edge [bend right=30,blue] node {} (5)
edge [bend right = 10,blue] node {} (5)
edge [bend left=30, blue] node {} (5)
edge [bend left = 10, blue] node{} (5)
edge node{} (1)

(3) edge [bend right=30, blue] node {} (6)
edge [bend right = 10,blue] node {} (6)
edge [bend left=30, blue] node {} (6)
edge [bend left = 10,blue] node{} (6)
edge [yellow] node {} (4)

(4) edge [bend right=30, blue] node {} (7)
edge [bend right = 10,blue] node {} (7)
edge [bend left=30, blue] node {} (7)
edge [bend left = 10, blue] node{} (7)
edge [green] node {}(6)

(5) edge [bend right=30, blue] node {} (8)
edge [bend right = 10,blue] node {} (8)
edge [bend left=30, blue] node {} (8)
edge [bend left =10, blue] node{} (8)
edge [red] node {}(3)
edge node{} (4)

(6) edge [bend right, blue] node {} (9)
edge [bend right =10,blue] node {} (9)
edge [bend left=30, blue] node {} (9)
edge [bend left=10,blue] node{} (9)
edge [yellow] node{} (7)

(7) edge [green] node {} (9)
edge [bend right=30, blue] node {} (10)
edge [bend right=10,blue] node {} (10)
edge [bend left=30, blue] node {} (10)
edge [bend left=10,blue] node{} (10)

(8) edge [red] node {}(6)
edge [bend right=30, blue] node {} (11)
edge [bend left=30,blue] node {} (11)
edge [bend left=10, blue] node {} (11)
edge [bend right=10,blue] node{} (11)
edge node{} (7)

(9) edge [bend right, blue] node {} (12)
edge [bend right =10,blue] node {} (12)
edge [bend left=30, blue] node {} (12)
edge [bend left=10,blue] node{} (12)
edge [yellow] node{} (10)

(10) edge [green] node {} (12)
edge [bend right=30, blue] node {} (13)
edge [bend right=10,blue] node {} (13)
edge [bend left=30, blue] node {} (13)
edge [bend left=10,blue] node{} (13)

(11) edge [red] node {}(9)
edge [bend right=30, blue] node {} (14)
edge [bend left=30,blue] node {} (14)
edge [bend left=10, blue] node {} (14)
edge [bend right=10,blue] node{} (14)
edge node{} (10)
(12) edge [yellow] node{} (13)

(14) edge node{} (13)
edge [red] node {} (12);

\end{tikzpicture}
$$ For a divisor $D \in Div_T(X) $ as above and $u \in \mathbb{C}$ consider the loop $$ \gamma^u_{D}(\theta) := \sum_{i=1}^{n-1} e^{-u+2 \pi i n'_i \theta } z_i + e^{2 \pi i m' \theta} w + \frac{e^{2 \pi i n'' \theta }}{ \prod_{i=1}^{n-1} z_i \cdot w^b} + \frac{e^{2 \pi i m'' \theta+ \pi i} }{\prod_{i=1}^{n-1} z_i \cdot w^{b+1} } + \frac{e^{2 \pi i l'' \theta}}{w} $$ For $ \theta \in [0,1)$. Define $ \Gamma_D : = [\gamma_D^t]$ for $0<<t$ and set $\widetilde{M}_D := M(\Gamma_D) \in Aut(Crit(X_{n,b} ; f_t))$. Let $$Crit(X) = \left \{ (z_{k},w_{l} ) \right \}_{k=0,l=0}^{n-1,1} \bigcup \left \{ (z'_m,w'_m) \right \}_{m=1}^{n-1} \subset (\mathbb{C}^{\ast})^{n-1} \times \mathbb{C}^{\ast} $$ be as in Proposition 3.3.1. The following example illustrates the $M$-aligned property for the class (c):
\bigskip

\hspace{-0.6cm} \bf Example 4.6 \rm (monodromies for $X= Bl_B(\mathcal{O}_{\mathbb{P}^4} \oplus \mathcal{O}_{\mathbb{P}^4} )$ with $B= \mathbb{P}^3$): The exceptional map of Definition 3.3.2 is defined in terms of the solutions of the LG-system for $$ f_t(z,w) = e^{-t} z_1+e^{-t}z_2+e^{-t}z_3+e^{-t}z_4+e^{-t} w + \frac{1}{z_1z_2 z_3 z_4} + \frac{i}{z_1 z_2 z_3 z_4 w } + \frac{1}{w} $$ for $0<<t$. The following is a graphical description of the $z_i$-coordinate of the solution scheme $Crit(X ; f_t)$:

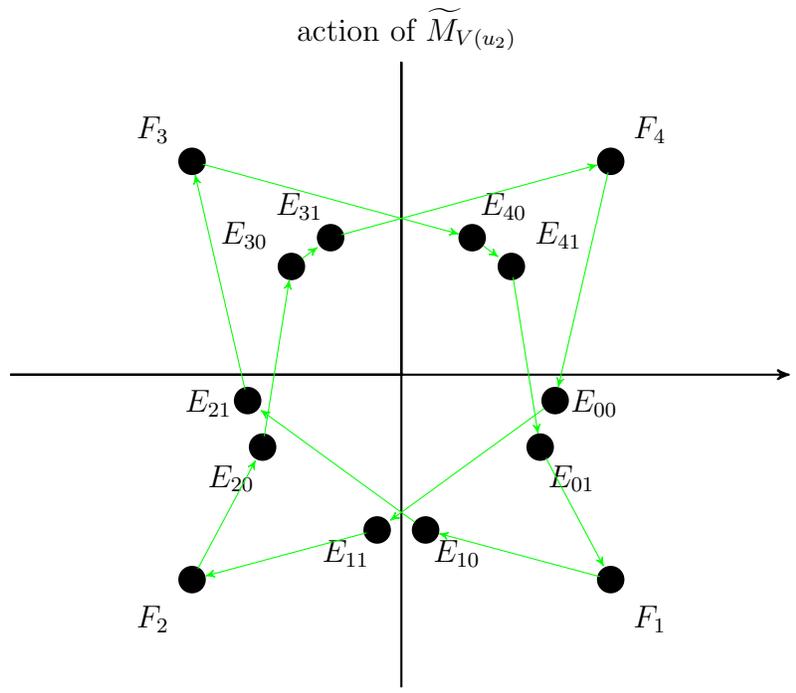
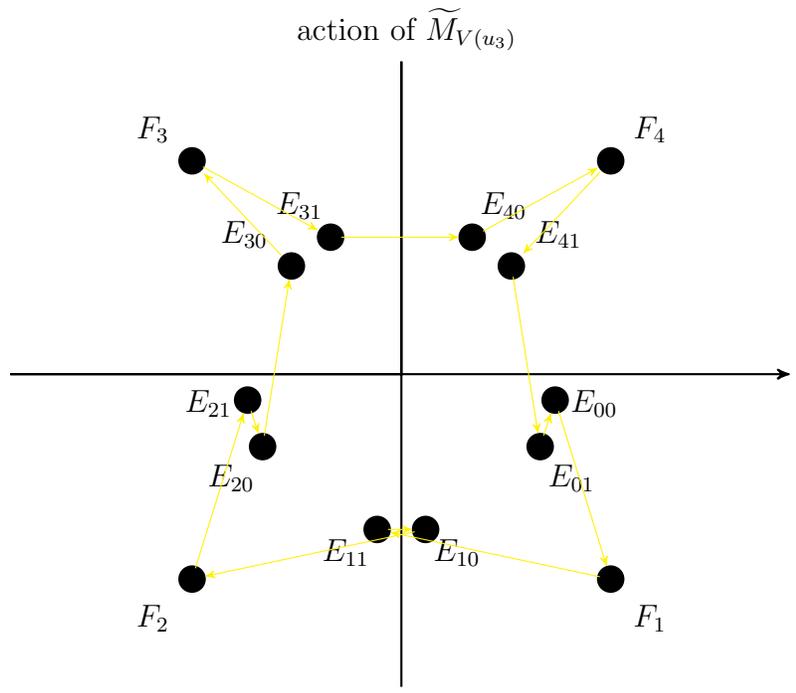
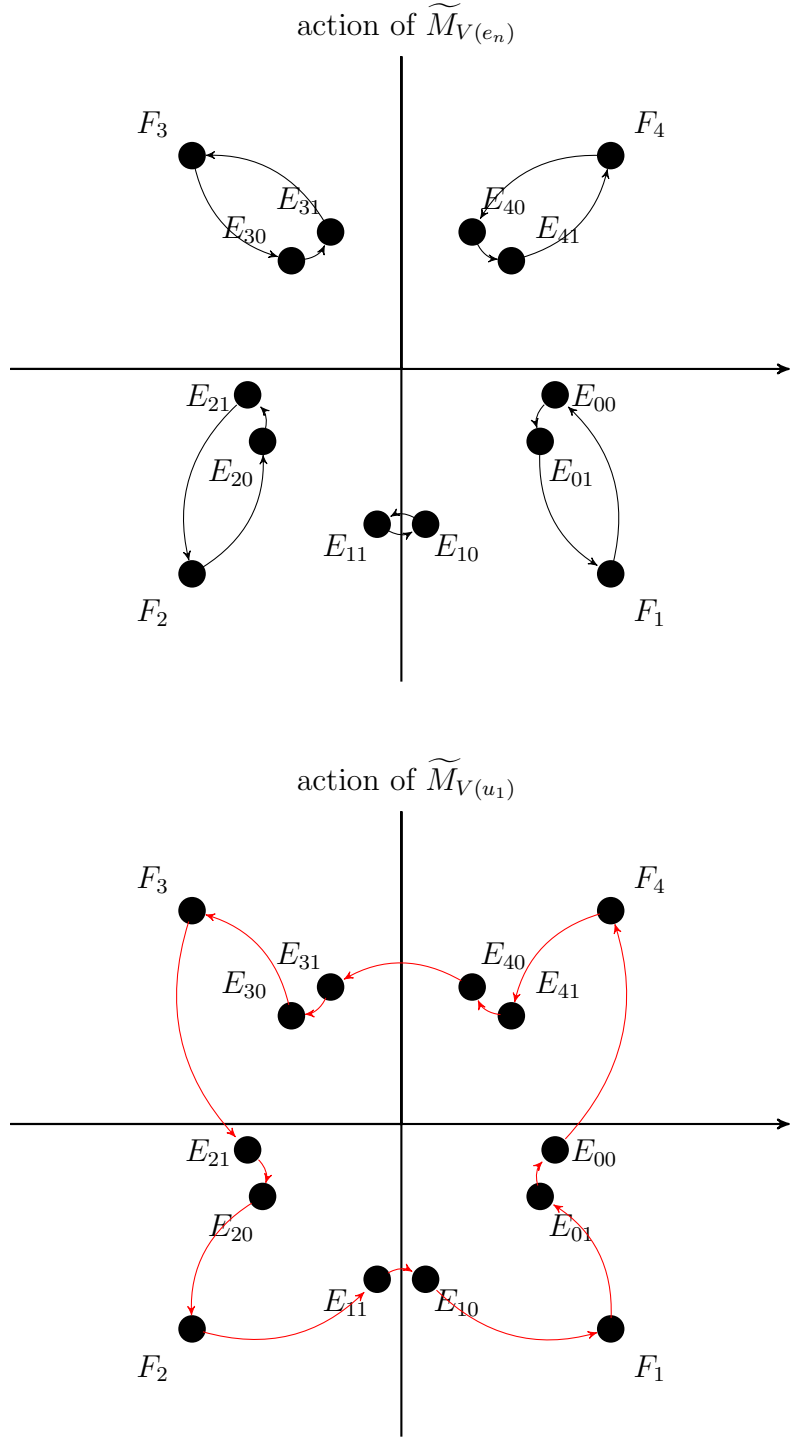
$$ \begin{tikzpicture}[->,>=stealth',shorten >=1pt,auto,node distance=3.5cm,main node/.style={circle,draw,font=\sffamily\bfseries \small}, scale=10.4]
\draw [->,thick] (0,-0.4)--(0,0.4) node (yaxis) [above] {$z(Crit(X; f_t))$}
|- (-0.5,0)--(0.5,0) node (xaxis) [right] {};

\node (a) at (0.19662, -0.0331365){};
\node (b) at (-0.19662, -0.0331365){};
\node (c) at (0.0309771, -0.198729){};
\node (d) at (-0.0309771, -0.198729){};
\node (e) at (0.0905276, 0.175303){};
\node (f) at (-0.0905276, 0.175303){};
\node (g) at (0.267681, -0.262253){};
\node (h) at (-0.267681, -0.262253){};
\node (h1) at (0.267703, 0.273023){};
\node (h2) at (-0.267703, 0.273023){};
\node (h3) at (0.177486, -0.0926983){};
\node (h4) at (-0.177486, -0.0926983){};
\node (h5) at (0.140639, 0.13849){};
\node (h6) at (-0.140639, 0.13849){};

\fill[black] (a) circle (0.5pt);
\fill[black] (b) circle (0.5pt);
\fill[black] (c) circle (0.5pt);
\fill[black] (d) circle (0.5pt);
\fill[black] (e) circle (0.5pt);
\fill[black] (f) circle (0.5pt);
\fill[black] (g) circle (0.5pt);
\fill[black] (h) circle (0.5pt);
\fill[black] (h1) circle (0.5pt);
\fill[black] (h2) circle (0.5pt);
\fill[black] (h3) circle (0.5pt);
\fill[black] (h4) circle (0.5pt);
\fill[black] (h5) circle (0.5pt);
\fill[black] (h6) circle (0.5pt);

\node (a1) at (0.24662, -0.0371365){$E_{00}$};
\node (b1) at (-0.24662, -0.0371365){$E_{21}$};
\node (c1) at (0.0709771, -0.228729){$E_{10}$};
\node (d1) at (-0.0709771, -0.228729){$E_{11}$};
\node (e1) at (0.1305276, 0.215303){$E_{40}$};
\node (f1) at (-0.1305276, 0.215303){$E_{31}$};
\node (g1) at (0.317681, -0.312253){$F_1$};
\node (h01) at (-0.317681, -0.312253){$F_2$};
\node (h11) at (0.317703, 0.313023){$F_4$};
\node (h12) at (-0.317703, 0.313023){$F_3$};
\node (h13) at (0.217486, -0.1326983){$E_{01}$};
\node (h14) at (-0.217486, -0.1326983){$E_{20}$};
\node (h15) at (0.200639, 0.17849){$E_{41}$};
\node (h16) at (-0.200639, 0.17849){$E_{30}$ };

\end{tikzpicture}$$

\hspace{-0.6cm} The corresponding monodromies are described as follows:

$$ \begin{tikzpicture}[->,>=stealth',shorten >=1pt,auto,node distance=3.5cm,main node/.style={circle,draw,font=\sffamily\bfseries \small}, scale=10.4]
\draw [->,thick] (0,-0.4)--(0,0.4) node (yaxis) [above] {$\textrm{ action of } \widetilde{M}_{V(e_i)} \textrm{ for } i=1,...,n-1$}
|- (-0.5,0)--(0.5,0) node (xaxis) [right] {};

\node (a) at (0.19662, -0.0331365){};
\node (b) at (-0.19662, -0.0331365){};
\node (c) at (0.0309771, -0.198729){};
\node (d) at (-0.0309771, -0.198729){};
\node (e) at (0.0905276, 0.175303){};
\node (f) at (-0.0905276, 0.175303){};
\node (g) at (0.267681, -0.262253){};
\node (h) at (-0.267681, -0.262253){};
\node (h1) at (0.267703, 0.273023){};
\node (h2) at (-0.267703, 0.273023){};
\node (h3) at (0.177486, -0.0926983){};
\node (h4) at (-0.177486, -0.0926983){};
\node (h5) at (0.140639, 0.13849){};
\node (h6) at (-0.140639, 0.13849){};

\fill[black] (a) circle (0.5pt);
\fill[black] (b) circle (0.5pt);
\fill[black] (c) circle (0.5pt);
\fill[black] (d) circle (0.5pt);
\fill[black] (e) circle (0.5pt);
\fill[black] (f) circle (0.5pt);
\fill[black] (g) circle (0.5pt);
\fill[black] (h) circle (0.5pt);
\fill[black] (h1) circle (0.5pt);
\fill[black] (h2) circle (0.5pt);
\fill[black] (h3) circle (0.5pt);
\fill[black] (h4) circle (0.5pt);
\fill[black] (h5) circle (0.5pt);
\fill[black] (h6) circle (0.5pt);

\node (a1) at (0.24662, -0.0371365){$E_{00}$};
\node (b1) at (-0.24662, -0.0371365){$E_{21}$};
\node (c1) at (0.0709771, -0.228729){$E_{10}$};
\node (d1) at (-0.0709771, -0.228729){$E_{11}$};
\node (e1) at (0.1305276, 0.215303){$E_{40}$};
\node (f1) at (-0.1305276, 0.215303){$E_{31}$};
\node (g1) at (0.317681, -0.312253){$F_1$};
\node (h01) at (-0.317681, -0.312253){$F_2$};
\node (h11) at (0.317703, 0.313023){$F_4$};
\node (h12) at (-0.317703, 0.313023){$F_3$};
\node (h13) at (0.217486, -0.1326983){$E_{01}$};
\node (h14) at (-0.217486, -0.1326983){$E_{20}$};
\node (h15) at (0.200639, 0.17849){$E_{41}$};
\node (h16) at (-0.200639, 0.17849){$E_{30}$ };

\path[every node/.style={font=\sffamily\small}]
(h) edge [blue, ->,bend left] (h2)
(h2) edge [blue,->,bend left] (h1)
(h1) edge [blue,->,bend left] (g)
(g) edge [blue,->,bend left] (h)

(b) [blue] edge (f)
(f) [blue] edge (h5)
(h5) [blue] edge (h3)
(h3) [blue] edge (d)
(d) [blue] edge (b)

(h4) [blue] edge (h6)
(h6) [blue] edge (e)
(e) [blue] edge (a)
(a) [blue] edge (c)
(c) [blue] edge (h4);

\end{tikzpicture}$$

$$ \begin{tikzpicture}[->,>=stealth',shorten >=1pt,auto,node distance=3.5cm,main node/.style={circle,draw,font=\sffamily\bfseries \small}, scale=10.4]
\draw [->,thick] (0,-0.4)--(0,0.4) node (yaxis) [above] {$\textrm{ action of } \widetilde{M}_{V(u_2)}$}
|- (-0.5,0)--(0.5,0) node (xaxis) [right] {};

\node (a) at (0.19662, -0.0331365){};
\node (b) at (-0.19662, -0.0331365){};
\node (c) at (0.0309771, -0.198729){};
\node (d) at (-0.0309771, -0.198729){};
\node (e) at (0.0905276, 0.175303){};
\node (f) at (-0.0905276, 0.175303){};
\node (g) at (0.267681, -0.262253){};
\node (h) at (-0.267681, -0.262253){};
\node (h1) at (0.267703, 0.273023){};
\node (h2) at (-0.267703, 0.273023){};
\node (h3) at (0.177486, -0.0926983){};
\node (h4) at (-0.177486, -0.0926983){};
\node (h5) at (0.140639, 0.13849){};
\node (h6) at (-0.140639, 0.13849){};

\fill[black] (a) circle (0.5pt);
\fill[black] (b) circle (0.5pt);
\fill[black] (c) circle (0.5pt);
\fill[black] (d) circle (0.5pt);
\fill[black] (e) circle (0.5pt);
\fill[black] (f) circle (0.5pt);
\fill[black] (g) circle (0.5pt);
\fill[black] (h) circle (0.5pt);
\fill[black] (h1) circle (0.5pt);
\fill[black] (h2) circle (0.5pt);
\fill[black] (h3) circle (0.5pt);
\fill[black] (h4) circle (0.5pt);
\fill[black] (h5) circle (0.5pt);
\fill[black] (h6) circle (0.5pt);

\node (a1) at (0.24662, -0.0371365){$E_{00}$};
\node (b1) at (-0.24662, -0.0371365){$E_{21}$};
\node (c1) at (0.0709771, -0.228729){$E_{10}$};
\node (d1) at (-0.0709771, -0.228729){$E_{11}$};
\node (e1) at (0.1305276, 0.215303){$E_{40}$};
\node (f1) at (-0.1305276, 0.215303){$E_{31}$};
\node (g1) at (0.317681, -0.312253){$F_1$};
\node (h01) at (-0.317681, -0.312253){$F_2$};
\node (h11) at (0.317703, 0.313023){$F_4$};
\node (h12) at (-0.317703, 0.313023){$F_3$};
\node (h13) at (0.217486, -0.1326983){$E_{01}$};
\node (h14) at (-0.217486, -0.1326983){$E_{20}$};
\node (h15) at (0.200639, 0.17849){$E_{41}$};
\node (h16) at (-0.200639, 0.17849){$E_{30}$ };

\path[every node/.style={font=\sffamily\small}]
(h2) edge [green] (e)
(f) edge [green] (h1)
(e) edge [green] (h5)
(h5) edge [green] (h3)
(a) edge [green] (d)
(h3) edge [green] (g)
(g) edge [green] (c)
(d) edge [green] (h)
(c) edge [green] (b)
(h) edge [green] (h4)
(h4) edge [green] (h6)
(h6) edge [green] (f)
(b) edge [green] (h2)
(h1) edge [green] (a);
\end{tikzpicture}$$

$$ \begin{tikzpicture}[->,>=stealth',shorten >=1pt,auto,node distance=3.5cm,main node/.style={circle,draw,font=\sffamily\bfseries \small}, scale=10.4]
\draw [->,thick] (0,-0.4)--(0,0.4) node (yaxis) [above] {$\textrm{ action of } \widetilde{M}_{V(u_3)}$}
|- (-0.5,0)--(0.5,0) node (xaxis) [right] {};

\node (a) at (0.19662, -0.0331365){};
\node (b) at (-0.19662, -0.0331365){};
\node (c) at (0.0309771, -0.198729){};
\node (d) at (-0.0309771, -0.198729){};
\node (e) at (0.0905276, 0.175303){};
\node (f) at (-0.0905276, 0.175303){};
\node (g) at (0.267681, -0.262253){};
\node (h) at (-0.267681, -0.262253){};
\node (h1) at (0.267703, 0.273023){};
\node (h2) at (-0.267703, 0.273023){};
\node (h3) at (0.177486, -0.0926983){};
\node (h4) at (-0.177486, -0.0926983){};
\node (h5) at (0.140639, 0.13849){};
\node (h6) at (-0.140639, 0.13849){};

\fill[black] (a) circle (0.5pt);
\fill[black] (b) circle (0.5pt);
\fill[black] (c) circle (0.5pt);
\fill[black] (d) circle (0.5pt);
\fill[black] (e) circle (0.5pt);
\fill[black] (f) circle (0.5pt);
\fill[black] (g) circle (0.5pt);
\fill[black] (h) circle (0.5pt);
\fill[black] (h1) circle (0.5pt);
\fill[black] (h2) circle (0.5pt);
\fill[black] (h3) circle (0.5pt);
\fill[black] (h4) circle (0.5pt);
\fill[black] (h5) circle (0.5pt);
\fill[black] (h6) circle (0.5pt);

\node (a1) at (0.24662, -0.0371365){$E_{00}$};
\node (b1) at (-0.24662, -0.0371365){$E_{21}$};
\node (c1) at (0.0709771, -0.228729){$E_{10}$};
\node (d1) at (-0.0709771, -0.228729){$E_{11}$};
\node (e1) at (0.1305276, 0.215303){$E_{40}$};
\node (f1) at (-0.1305276, 0.215303){$E_{31}$};
\node (g1) at (0.317681, -0.312253){$F_1$};
\node (h01) at (-0.317681, -0.312253){$F_2$};
\node (h11) at (0.317703, 0.313023){$F_4$};
\node (h12) at (-0.317703, 0.313023){$F_3$};
\node (h13) at (0.217486, -0.1326983){$E_{01}$};
\node (h14) at (-0.217486, -0.1326983){$E_{20}$};
\node (h15) at (0.200639, 0.17849){$E_{41}$};
\node (h16) at (-0.200639, 0.17849){$E_{30}$ };

\path[every node/.style={font=\sffamily\small}]
(h2) edge [yellow] (f)
(f) edge [yellow] (e)
(e) edge [yellow] (h1)
(h1) edge [yellow] (h5)
(h5) edge [yellow] (h3)
(h3) edge [yellow] (a)
(a) edge [yellow] (g)
(g) edge [yellow] (d)
(d) edge [yellow] (c)
(c) edge [yellow] (h)
(h) edge [yellow] (b)
(b) edge [yellow] (h4)
(h4) edge [yellow] (h6)
(h6) edge [yellow] (h2);

\end{tikzpicture}$$

$$ \begin{tikzpicture}[->,>=stealth',shorten >=1pt,auto,node distance=3.5cm,main node/.style={circle,draw,font=\sffamily\bfseries \small}, scale=10.4]
\draw [->,thick] (0,-0.4)--(0,0.4) node (yaxis) [above] {$\textrm{ action of } \widetilde{M}_{V(e_n)}$}
|- (-0.5,0)--(0.5,0) node (xaxis) [right] {};

\node (a) at (0.19662, -0.0331365){};
\node (b) at (-0.19662, -0.0331365){};
\node (c) at (0.0309771, -0.198729){};
\node (d) at (-0.0309771, -0.198729){};
\node (e) at (0.0905276, 0.175303){};
\node (f) at (-0.0905276, 0.175303){};
\node (g) at (0.267681, -0.262253){};
\node (h) at (-0.267681, -0.262253){};
\node (h1) at (0.267703, 0.273023){};
\node (h2) at (-0.267703, 0.273023){};
\node (h3) at (0.177486, -0.0926983){};
\node (h4) at (-0.177486, -0.0926983){};
\node (h5) at (0.140639, 0.13849){};
\node (h6) at (-0.140639, 0.13849){};

\fill[black] (a) circle (0.5pt);
\fill[black] (b) circle (0.5pt);
\fill[black] (c) circle (0.5pt);
\fill[black] (d) circle (0.5pt);
\fill[black] (e) circle (0.5pt);
\fill[black] (f) circle (0.5pt);
\fill[black] (g) circle (0.5pt);
\fill[black] (h) circle (0.5pt);
\fill[black] (h1) circle (0.5pt);
\fill[black] (h2) circle (0.5pt);
\fill[black] (h3) circle (0.5pt);
\fill[black] (h4) circle (0.5pt);
\fill[black] (h5) circle (0.5pt);
\fill[black] (h6) circle (0.5pt);

\node (a1) at (0.24662, -0.0371365){$E_{00}$};
\node (b1) at (-0.24662, -0.0371365){$E_{21}$};
\node (c1) at (0.0709771, -0.228729){$E_{10}$};
\node (d1) at (-0.0709771, -0.228729){$E_{11}$};
\node (e1) at (0.1305276, 0.215303){$E_{40}$};
\node (f1) at (-0.1305276, 0.215303){$E_{31}$};
\node (g1) at (0.317681, -0.312253){$F_1$};
\node (h01) at (-0.317681, -0.312253){$F_2$};
\node (h11) at (0.317703, 0.313023){$F_4$};
\node (h12) at (-0.317703, 0.313023){$F_3$};
\node (h13) at (0.217486, -0.1326983){$E_{01}$};
\node (h14) at (-0.217486, -0.1326983){$E_{20}$};
\node (h15) at (0.200639, 0.17849){$E_{41}$};
\node (h16) at (-0.200639, 0.17849){$E_{30}$ };

\path[every node/.style={font=\sffamily\small}]
(h1) edge [ ->,bend right] (e)
(e) edge [->,bend right] (h5)
(h5) edge [->,bend right] (h1)

(g) edge [ ->,bend right] (a)
(a) edge [->,bend right] (h3)
(h3) edge [->,bend right] (g)

(b) edge [ ->,bend right] (h)
(h) edge [->,bend right] (h4)
(h4) edge [->,bend right] (b)

(h2) edge [ ->,bend right] (h6)
(h6) edge [->,bend right] (f)
(f) edge [->,bend right] (h2)

(c) edge [ ->,bend right] (d)
(d) edge [->,bend right] (c);

\end{tikzpicture}$$

$$ \begin{tikzpicture}[->,>=stealth',shorten >=1pt,auto,node distance=3.5cm,main node/.style={circle,draw,font=\sffamily\bfseries \small}, scale=10.4]
\draw [->,thick] (0,-0.4)--(0,0.4) node (yaxis) [above] {$\textrm{ action of } \widetilde{M}_{V(u_1)}$}
|- (-0.5,0)--(0.5,0) node (xaxis) [right] {};

\node (a) at (0.19662, -0.0331365){};
\node (b) at (-0.19662, -0.0331365){};
\node (c) at (0.0309771, -0.198729){};
\node (d) at (-0.0309771, -0.198729){};
\node (e) at (0.0905276, 0.175303){};
\node (f) at (-0.0905276, 0.175303){};
\node (g) at (0.267681, -0.262253){};
\node (h) at (-0.267681, -0.262253){};
\node (h1) at (0.267703, 0.273023){};
\node (h2) at (-0.267703, 0.273023){};
\node (h3) at (0.177486, -0.0926983){};
\node (h4) at (-0.177486, -0.0926983){};
\node (h5) at (0.140639, 0.13849){};
\node (h6) at (-0.140639, 0.13849){};

\fill[black] (a) circle (0.5pt);
\fill[black] (b) circle (0.5pt);
\fill[black] (c) circle (0.5pt);
\fill[black] (d) circle (0.5pt);
\fill[black] (e) circle (0.5pt);
\fill[black] (f) circle (0.5pt);
\fill[black] (g) circle (0.5pt);
\fill[black] (h) circle (0.5pt);
\fill[black] (h1) circle (0.5pt);
\fill[black] (h2) circle (0.5pt);
\fill[black] (h3) circle (0.5pt);
\fill[black] (h4) circle (0.5pt);
\fill[black] (h5) circle (0.5pt);
\fill[black] (h6) circle (0.5pt);

\node (a1) at (0.24662, -0.0371365){$E_{00}$};
\node (b1) at (-0.24662, -0.0371365){$E_{21}$};
\node (c1) at (0.0709771, -0.228729){$E_{10}$};
\node (d1) at (-0.0709771, -0.228729){$E_{11}$};
\node (e1) at (0.1305276, 0.215303){$E_{40}$};
\node (f1) at (-0.1305276, 0.215303){$E_{31}$};
\node (g1) at (0.317681, -0.312253){$F_1$};
\node (h01) at (-0.317681, -0.312253){$F_2$};
\node (h11) at (0.317703, 0.313023){$F_4$};
\node (h12) at (-0.317703, 0.313023){$F_3$};
\node (h13) at (0.217486, -0.1326983){$E_{01}$};
\node (h14) at (-0.217486, -0.1326983){$E_{20}$};
\node (h15) at (0.200639, 0.17849){$E_{41}$};
\node (h16) at (-0.200639, 0.17849){$E_{30}$ };

\path[every node/.style={font=\sffamily\small}]
(a) edge [red, bend right] (h1)
(h1) edge [red, bend right] (h5)
(h5) edge [red, bend left] (e)
(e) edge [red, bend right] (f)
(f) edge [red, bend left] (h6)
(h6) edge [red, bend right] (h2)
(h2) edge [red, bend right] (b)
(b) edge [red, bend left] (h4)
(h4) edge [red, bend right] (h)
(h) edge [red, bend right] (d)
(d) edge [red, bend left] (c)
(c) edge [red, bend right] (g)
(g) edge [red, bend right] (h3)
(h3) edge [red, bend left] (a);

\end{tikzpicture}$$ In particular, one has:

\bigskip

(1) $\widetilde{M}_{V(e_i)}(w_k,w_l)=(z_{k+1},w_l)$ and $ \widetilde{M}_{V(e_i)}(z'_m,w'_m) = (z'_{m+1},w'_m)$ for $i=0,...,n-1$.

\bigskip

(2) $\widetilde{M}_{V(e_n)}(z_k,w_0)=(z_{k-b},w_{1})$ for $b+1 \leq k \leq n-1$.
\bigskip

(3) $\widetilde{M}_{V(u_1)}(z_k,w_0)=(z'_{k-b},w'_k)$ for $b+1 \leq k \leq n-1$.

\bigskip

(4) $\widetilde{M}_{V(u_2)}(z_k,w_1)=(z'_{k+1},w'_{k+1})$ for $0 \leq k \leq n-2$.
\bigskip

(5) $\widetilde{M}_{V(u_3)}(z'_k,w'_k)=(z_k,w_1)$ for $1 \leq k \leq n-1$.

\bigskip

\hspace{-0.6cm} Relations (1)-(5) establish the $M$-aligned property in this case.

\section{Discussion and Concluding Remarks}
\label{s:LGM}

\bigskip

\hspace{-0.6cm} We would like to conclude with the following remarks and questions:

\bigskip

\hspace{-0.6cm} \bf (a) Monodromies and Lagrangian submanifolds: \rm A leading source of interest for the study of the structure of $\mathcal{D}^b(X)$, in recent years, has been their role in the famous homological mirror symmetry conjecture due to Kontsevich, see \cite{K}. For a toric Fano manifold, the analog statement of homological mirror symmetry relates the structure of $\mathcal{D}^b(X)$ to the structure of the Fukaya-Seidel cateogry $Fuk(\widetilde{Y}^{\circ})$, where $\widetilde{Y}^{\circ}$ is a disingularization of a hyperplane section $Y^{\circ}$ of $X^{\circ}$, see \cite{AKO,Ko2,S}. In this setting, homological mirror symmetry was proven by Abouzaid in \cite{Ab1,Ab2}. On the other hand, the Fukaya-Seidel category is known to admit exceptional collections of Lagrangian spheres, due Seidel's Lefschetz thimble construction, see for instance \cite{S2}. Thus, via homological mirror symmetry these collections translate to exceptional collections of objects in $\mathcal{D}^b(X)$. However, it seems a non-trivial question to discern, in terms of the Fukay-Seidel category itself, which such collections of vanishing spheres (if any) in $Fuk(\widetilde{Y}^{\circ})$ actually correspond to collections of line bundles in $\mathcal{D}^b(X)$, under the mirror map. Examples could be found in \cite{AKO}. It is thus interesting to pose the following question:

\bigskip

\hspace{-0.6cm} \bf Question: \rm Is it possible to naturally associate Lagrangian submanifolds $L(z) \subset \widetilde{Y}^{\circ}$ to solutions of the Landau-Ginzburg system $z \in Crit(X)$ with the property $$HF(L(z),L(w)) \simeq Hom_{mon}(z,w) \hspace{0.5cm} \textrm{ for} \hspace{0.25cm} z,w \in Crit(X)$$ where $HF$ stands for Lagrangian Floer homology?

\bigskip

\hspace{-0.6cm} (b) \bf Further toric Fano manifolds: \rm It is interesting to ask to which extent the methods of this work could be generalized to further classes of torc Fano manifolds? In particular, the main question seems to be in which cases there exists an asymptotic deformation $f_t$ of the Landau-Ginzburg potential $f_X$ such that the corresponding map $E : Crit(X ; f_t) \rightarrow Pic(X)$, defined in section 3, is exceptional for $0<<t$.

\hspace{-0.6cm} Regrading the deformation procedure, let us note that the zero set $V(f)= \left \{f=0 \right \} \subset (\mathbb{C}^{\ast})^n$ of an element $f \in L(\Delta^{\circ})$
is an affine Calabi-Yau hypersurface. In \cite{Ba2} Batyrev introduced $\mathcal{M}(\Delta^{\circ})$ the toric moduli of such affine Calabi-Yau hyper-surfaces which is a $\rho(X)$-dimensional singular toric variety obtained as the quotient of $L(\Delta^{\circ})$ by appropriate equivalence relations. In \cite{Ba2} Batyrev further shows that $PH^{n-1}(V(f)) \simeq Jac(f)$, where $Jac(f)$ is the function ring of the solution scheme $Crit(X ; f) \subset (\mathbb{C}^{\ast})^n$.

\hspace{-0.6cm} In this sense our approach could be viewed as a suggesting that in the toric Fano case homological data about the structure of $\mathcal{D}^b(X)$, could, in fact, be extracted from the local behavior around the boundary of the B-model moduli, which in our case is $\mathcal{M}(\Delta^{\circ})$, rather than the Fukaya category appearing in the general homological mirror symmetry conjecture, whose structure is typically much harder to analyze.

\bigskip

\hspace{-0.6cm} \bf Acknowledgements: \rm The author would like to thank Leonid Polterovich for his ongoing support and interest in this project. The author would also like to thank Denis Auroux for valuable comments and suggestions on previous versions of this work. This research has been partially supported by the European
Research Council Advanced grant 338809 and AdG 669655.

\end{document}